\documentclass[12pt]{article}

\usepackage{color}
\usepackage{amsmath}
\usepackage{amssymb}
\usepackage{graphicx}
\usepackage{float}
\usepackage{hyperref}
\usepackage{soul}

\def\Res{\mathop{\mathrm{Res}}}

\def\L{{\mathbf{L}}}
\def\bQ{{\mathbf{Q}}}

\def\cO{\mathcal{O}}

\def\qed{{\hfill $\Diamond$} \bigskip}

\def\rPic{\mathrm{Pic}}

\def\ch{{\rm ch}}
\def\td{{\rm td}} 
\def\Td{{\rm Td}}
\def\C{\mathbb{C}}
\def\R{\mathbb{R}}
\def\Z{\mathbb{Z}}
\def\Q{\mathbb{Q}}
\def\d{\partial}
\def\dbar{{\overline{\partial}}}

\def\isom{\xrightarrow{\sim}}
\def\Hom{{\rm Hom}}
\def\CP{\mathbb{C} \rm{P}}

\def\P{\mathbb{P}}

\usepackage{theorem}

\newtheorem{theorem}{Theorem}
\newtheorem{proposition}{Proposition}[section]
\newtheorem{corollary}[proposition]{Corollary}
\newtheorem{lemma}[proposition]{Lemma}

{\theorembodyfont{\rmfamily}

\newtheorem{example}[proposition]{Example}
\newtheorem{remark}[proposition]{Remark}
\newtheorem{notation}[proposition]{Notation}
}
\usepackage{soul}

\usepackage{lipsum}

\newcommand\blfootnote[1]{%
  \begingroup
  \renewcommand\thefootnote{}\footnote{#1}%
  \addtocounter{footnote}{-1}%
  \endgroup
}

\author{Semyon Klevtsov$^1$, Dimitri Zvonkine$^2$}

\title{The Chern character of the Laughlin vector bundle in the Fractional Quantum Hall Effect}

 \date{\today}

\begin{document}

\maketitle
{\it\small
\noindent
$^1$ IRMA, Université de Strasbourg, 7 r.~René Descartes, 67084 Strasbourg, France\\
$^2$ Université Paris-Saclay, UVSQ, CNRS, Laboratoire de Mathématiques de Versailles, 78000, Versailles, France.}
\blfootnote{The first author was partly supported by the the ANR-20-CE40-0017 grant, the Initiative
d’excellence program and the Institute for Advanced Study Fellowship of the University of Strasbourg.}

\begin{abstract}
We begin by explaining how a physical problem of studying the quantum Hall effect on a closed surface~$C$ leads, via Laughlin's approach, to a mathematical question of describing the rank and the first Chern class of a particular vector bundle on the Picard group $\rPic^g(C)$. Then we formulate and solve the problem mathematically, proving several important conjectures made by physicists, in particular the Wen-Niu topological degeneracy conjecture and the Wen-Zee shift formula.

Let $C$ be a closed Riemann surface of genus~$g$ and $S^NC$ its $N$th symmetric power. The product $C \times \rPic^d(C)$ carries a universal line bundle. On the product $C^N \times \rPic^d(C)$ we consider the product of $N$ pull-backs of this universal line bundle and twist it by a power of the diagonal on $C^N$. The resulting line bundle descends onto $S^NC \times \rPic^d(C)$. Its push-forward (as a sheaf) to $\rPic^d(C)$ is a vector bundle that we call Laughlin's vector bundle. We determine all the Chern characters of the Laughlin vector bundle via a Grothendieck-Riemann-Roch calculation. 
\end{abstract}

\section{Introduction: from the quantum Hall effect to algebraic geometry}
\label{Sec:PhysToMath}

%The target audience of this a section is a mathematician with some basic knowledge of quantum mechanics who would like to learn about the quantum Hall effect. 

\subsection{The Hall effect}
\label{Ssec:HallEffect}
The classical Hall effect was dicovered in the XIX century. Consider a thin conducting rectangular plate of size $L_x \times L_y$ with conductivity $\kappa$. Its full conductance along the directions $x$ and $y$ is then equal to $\frac{L_y}{L_x} \kappa$ and $\frac{L_x}{L_y} \kappa$ respectively. In other words, applying, for instance, a voltage $V_x$ between its left and right sides, we get a horizontal current $I_x = \frac{L_y}{L_x}\kappa V_x$. Thus the conductance matrix 
$$
\left(
\begin{array}{cc}
\frac{L_y}{L_x}\kappa & 0\\[10pt]
0 & \frac{L_x}{L_y}\kappa
\end{array}
\right)
$$
is diagonal.

Now apply a strong magnetic field perpendicular to the rectangle. The moving charged particles will be deviated by the Lorentz force, so a horizontal potential will lead to both a horizontal and a vertical current and the conductance matrix will have the form
$$
\left(
\begin{array}{cc}
\frac{L_y}{L_x}\kappa & \sigma_H \\[10pt]
-\sigma_H & \frac{L_x}{L_y}\kappa
\end{array}
\right)
$$
where $\sigma_H$ is called the {\em Hall conductance} and is independent of $L_x$ and $L_y$.

In certain materials at very low temperatures surprising quantum effects kick in: the Hall conductance $\sigma_H$, instead of growing linearly with $B$, presents more or less extended plateaux at values $\sigma_H = \nu\frac {e^2}h$, where the electronic charge squared divided by the Planck constant is a natural unit of conductance and $\nu$ is, to a very high precision, either an integer or a fractional number with a small denominator. As for the ordinary conductance $\kappa$, within the plateaux it vanishes, so that the conductance matrix becomes anti-diagonal. That is, a horizontal voltage creates a current only in the vertical direction and vice versa. 

The suggested explanation is that the ground state
%, {\color{red} called the lowest Landau level,} 
of the Hamiltonian of a charged particle on the conducting surface is highly degenerate (that is, the eigenspace with lowest eigenvalue has a high dimension). The quantum Hall effect occurs when it is completely filled with electrons, which must occupy different states since they are fermions.  Normally, this would require one to fix the intensity of the magnetic field very precisely, but due to impurities in the conductor, the effect is actually visible for a range of values of the magnetic field. For this reason, the quantum Hall effect is hard to observe in an extremely pure material.

(Alternatively, there can be several eigenstates of the Hamiltonian with energies very close to the ground state, separated by an energy gap from the other eigenstates. These eigenstates collectively from what is called the lowest Landau level.)

Simplest models of the quantum Hall effect identify opposite sides of the rectangle to create a torus:

\begin{figure}[H]
\begin{center}
\includegraphics[width=10em]{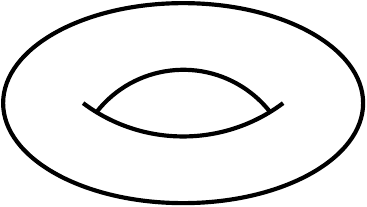} \hspace{2em}
\includegraphics[width=10em]{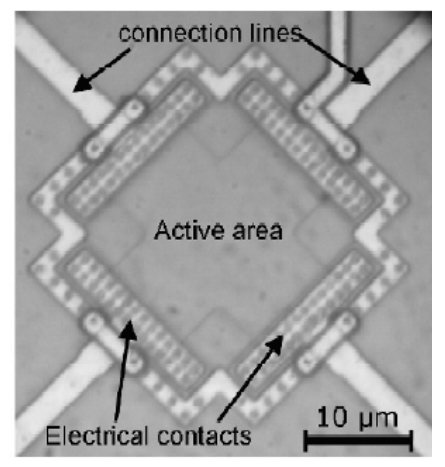}
\caption{Identifying opposite sides of a rectangle. Theory vs experiment (a Hall sensor from~\cite{MukHua}, p.~4).}
\end{center}
\end{figure}

In this paper, we will work in the following theoretical setup: our conductor is a closed surface of genus~$g$, and the number of particles is close to the number that is necessary to fill the ground state. In practice, a genus~2 experimental setup might look like this:

\begin{figure}[H]
\begin{center}
\includegraphics[width=20em]{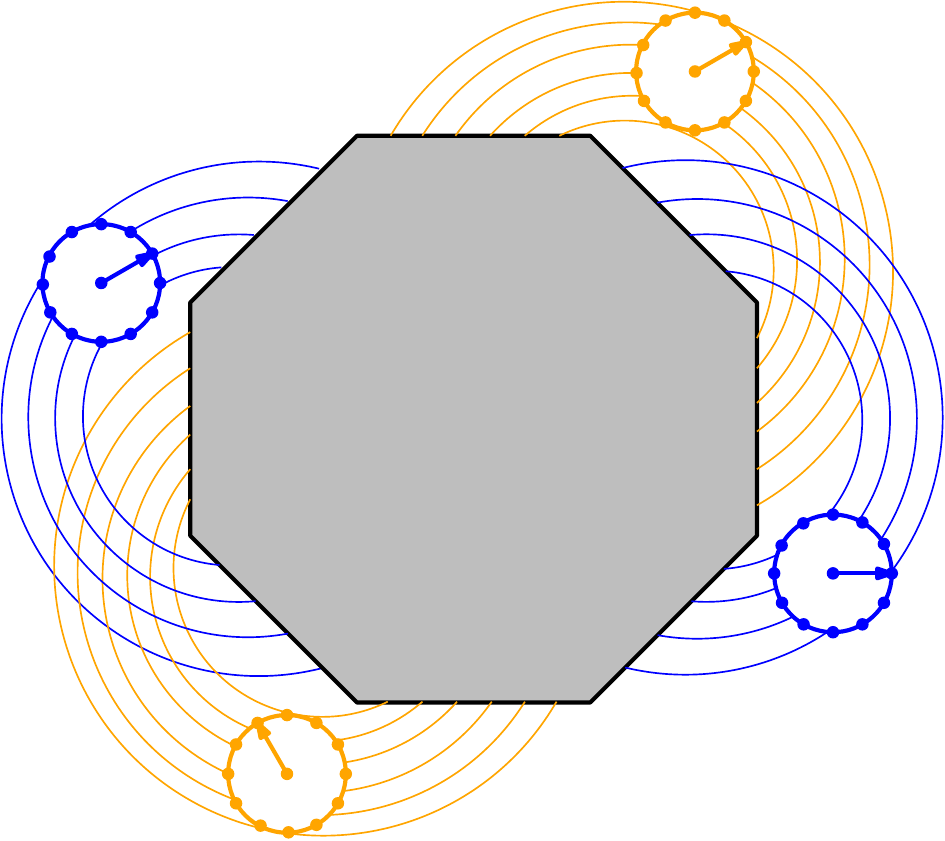}
\caption{Pairwise identification of sides of an octogon gives a genus~2 surface.}
\end{center}
\end{figure}

Note that the measuring devices are in a one-to-one correspondence with a basis of the cohomology group $H^1(C,\R)$.

\subsection{The setup}
\label{Ssec:Setup}

Let $C$ be a conducting surface with a Riemannian metric. In our models, $C$ will play the role of space, so that $C \times \R$ will be the space-time.
We represent the Riemannian metric as a pair $\omega$, $J$, where $\omega$ is a symplectic form and $J$ a complex structure. In particular, $(C,J)$ is a Riemann surface. 

A striking feature of this setup is that many physically meaningful effects only depend on the Riemann surface structure of~$C$, but not on~$\omega$. This allows one to study the quantum Hall effect using techniques from algebraic geometry. For instance, ground states of the Hamiltonian of one particle in the perpendicular magnetic field turn out to be holomorphic sections of a holomorphic line bundle on~$C$. Thus the degeneracy of the ground state is the dimension of the space of holomorphic sections of a line bundle, which is given by the Riemann-Roch formula. The main results of this paper can be viewed as a generalization of this observation, with $N$ particles instead of one, the electric field switched on, and the Grothendieck-Riemann-Roch formula replacing the Riemann-Roch formula.

An {\em electro-magnetic field} $(E,B)$ over $C$ is a principal $U(1)$-bundle $\mathbf{P} \to C \times \R$ over the space-time $C \times \R$, endowed with a connection~$\nabla$, up to a connection-preserving isomorphism (gauge transformation). 

Denote by $P \to C$ the restriction of the circle bundle $\mathbf{P}$ to $C \times \{0 \}$. The first Chern class of $P \to C$ is the total magnetic flux through~$C$. It is a topological invariant that cannot change with time. We will use the parallel transport along the time lines $z \times \R \subset C \times \R$ to identify $\mathbf{P}$ with $P \times \R$. Under this identification (Weyl gauge or temporal gauge in physical terminology), the connection $\nabla$ has the form $\nabla = d + \alpha_t$, where $\alpha_t$ is a time-dependent 1-form on~$C$. The connection has no $dt$ term due to the choice of trivialization. Thus $\nabla$ can be viewed as a family $\nabla_t = d + \alpha_t$ of connections on~$P$ depending on time. The electric and the magnetic fields are 
$$
E = \d \alpha_t / \d t,  \qquad B \omega = d \alpha_t.
$$
In other words, $E$ is a time-dependent 1-form on~$C$, and $B$ is a time-dependent function on~$C$ obtained by dividing the 2-form $d \alpha_t$ by the symplectic form~$\omega$.

From now on we will work under two assumptions. First, the magnetic field does not depend on~$t$. Second, there is no accumulation of charge on~$C$ and therefore the divergence of~$E$ vanishes. The first assumption is just the experimental setup. The second one is a good approximation for most conductors.

\subsection{A $2g$-dimensional parameter space for time}
\label{Ssec:2gtime}

The time-independent magnetic field assumption implies that the electric field $E_t$ is a closed 1-form on~$C$. Indeed $d \alpha_t = B \omega$ for all~$t$, therefore
$$
dE = d \frac{\d \alpha_t}{\d t} = 0.
$$
Since $E$ is, moreover, divergence-free, it is actually a harmonic 1-form. Note that the space of harmonic 1-forms on~$C$ only depends on the complex structure of~$C$: indeed, a harmonic 1-form is the real (or the imaginary) part of an abelian differential.

The space of harmonic 1-forms is naturally isomorphic to $H^1(C,\R)$; thus $H^1(C,\R)$ is the parameter space of all possible time evolutions of the connection $\nabla_t$. Moreover, we can construct a universal connection on the circle bundle $P \times H^1(C,\R)$. 

Let $P \to C$ be any degree~$d$ circle bundle with a connection $\nabla_0$. Then we define the connection $\nabla$ on the circle bundle $P \times H^1(C,\R) \to C \times H^1(C,\R)$ as
$$
\nabla(p, \alpha) = \nabla_0(p) + \alpha,
$$
where $\alpha$ in the left-hand side is a cohomology class, while in the right-hand side it is the corresponding harmonic 1-form. Every connection $\nabla_t$ on $\mathbf{P} \to C \times \R$ is the pullback of $\nabla$ under some map $\R \to H^1(C,\R)$. Note that the connection $\nabla$ has no component along $H^1(C,\R)$. This is compatible with our choice of gauge on $\mathbf{P} \to C \times \R$.

Actually, the universal connection $\nabla$ can be defined even more intrinsically without having to choose $\nabla_0$. Indeed, consider the space $W$ of connections on $P \to C$ up to gauge equivalence $\phi: C \to U(1)$ that can be lifted to $\widetilde{\phi}: C \to \R$. This is an affine space whose underlying vector space is $H^1(C,\R)$. There is a universal connection $\nabla$ on $P \times W$ such that (i)~its component along $W$ vanishes and (ii)~the difference between its restrictions to two fibers $P \times \{w\}$ and $P \times \{ w' \}$ is always a harmonic 1-form on~$C$. It is unique up to a choice of gauge transformation on any one fiber $P \times \{ w\}$.

\subsection{Factorizing by $H^1(C,\Z)$}
\label{Ssec:Factorc1}

The integer cohomology group $H^1(C,\Z)$ acts on the space $W$ by translations. Denote by $T = W/H^1(C,\Z)$ the quotient; it is a torus of real dimension~$2g$. For an $\alpha \in H^1(C,\Z)$ and a $w \in W$, the connections on $P \times \{w\}$ and $P \times \{ w + \alpha \}$ are gauge equivalent. But the gauge transformation $\phi = \exp\left(2 \pi i \int \!\!\alpha\right)$ cannot be lifted from $U(1)$ to $\R$ (unless $\alpha =0$).

Choose a base point $z_0 \in C$. Given this choice, we can lift the action of $H^1(C,\Z)$ to $P \times W$ via 
$$
\alpha: (z,w,\phi) \mapsto \left(z, w+\alpha, \phi \times \exp \left(2 \pi i \int_{z_0}^z \!\!\alpha \right) \right).
$$
Note, however, that this lifting depends on~$z_0$.

The quotient by the above action of $H^1(C,\Z)$ is a $U(1)$-bundle over $C \times T$ endowed with a connection. Changing the choice of~$z_0$ modifies the holonomy of the connection along a translation by $\alpha \in H^1(C,\Z)$ by a constant phase (that is, a phase independent of $z \in C$). There is no canonical way to settle this ambiguity. 

To sum up, we have obtained a circle bundle with a connection over $C \times T$, but the construction is not completely canonical: there are $2g$ holonomies that can be chosen arbitrarily. This choice will not affect the subsequent computations in any way.

\begin{figure}[H]
\begin{center}
\includegraphics[width=20em]{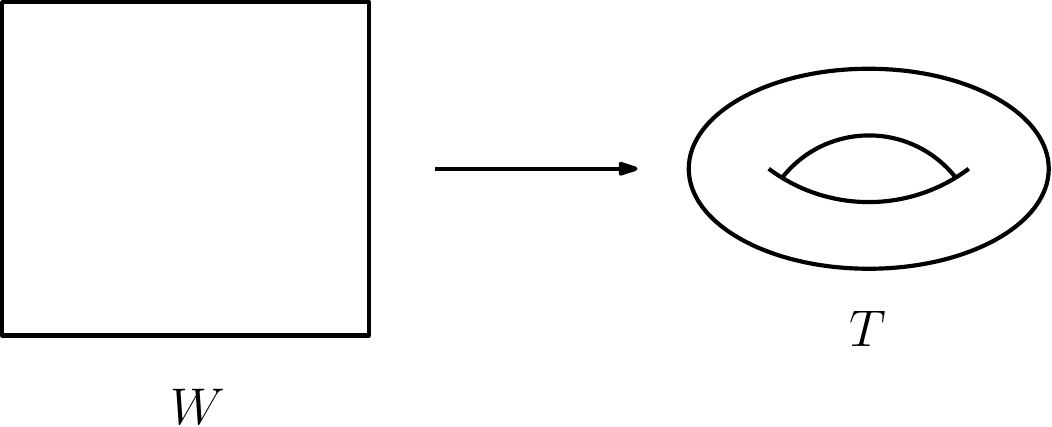}
\caption{A circle bundle with a canonical connection over $C \times W$ descends to a circle bundle over $C \times T$ with a connection defined up to a choice of phase $H^1(C,\Z) \to U(1)$.}
\end{center}
\end{figure}

To begin with, note that the first Chern class of 
$$
(P \times W) / H^1(C,\Z) \to C \times T
$$ 
does not depend on the choice of the holonomies. Let us write out this first Chern class. We have
$$
H^2(C \times T,\Z) = H^2(C,\Z)  + H^2(T,\Z) + H^1(C,\Z) \times H^1(T,\Z).
$$
Thus the first Chern class has three components. 

{\bf 1.} The first component can be read on the restriction $P \to C$ to one fiber. It equals $d \cdot \xi$, where $\xi \in H^2(C,\Z)$ is the Poincar\'e dual of the class of a point and, by abuse of notation, we also denote by $\xi$ its pull-back to $H^2(C \times T, \Z)$.

{\bf 2.} The second component vanishes, because the connection is trivial along~$T$.

{\bf 3.} To write out the third (mixed) component we need to introduce some notation. The group $H_1(T, \Z)$ is canonically isomorphic to $H^1(C,\Z)$. Thus there is a canonical map
$$
\begin{array}{rcccc}
\eta&:& H_1(T, \Z) \times H_1(C,\Z) &\to& \Z \\
&&\alpha,\gamma &\mapsto& (\alpha,\gamma)\\
\end{array}
$$
This map is an element of $H^1(T, \Z) \times H^1(C,\Z)$, and we claim that the mixed term of the first Chern class is the element of $H^2(T \times C,\Z)$ corresponding to~$\eta$. (Note that in the identification of $H^1(T, \Z) \times H^1(C,\Z)$ with $H^2(T \times C,\Z)$ the sign depends on the order of the factors.) This can be seen directly by applying the parallel transport along loops on $T \times C$. Let $\gamma$ be a 1-cycle on~$C$ and $\alpha$ a 1-form on~$C$ with integral periods corresponding to a 1-cycle on~$T$. The parallel transport along $\alpha$ from $w$ to $w + \alpha$ gives 0 holonomy. Then the parallel transport along~$\gamma$ over $w+\alpha \in W$ gives some holonomy that we denote by $\phi_{w+\alpha}(\gamma)$. Then the parallel transport from $w+\alpha$ to $w$ gives 0 holonomy. Finally, the parallel transport along $-\gamma$ over $w$ gives holonomy $-\phi_w(\gamma)$. The difference $\phi_{w+\alpha} - \phi_w$ equals the 1-form $\alpha$, so altogether we get the holonomy $(\alpha,\gamma)$. 

We have obtained the following statement.

\begin{proposition} \label{Prop:c1OneParticle}
The first Chern class of the $U(1)$-bundle
$$
(P \times W) / H^1(C,\Z) \to C \times T
$$ 
equals $d \cdot \xi + \eta$, where $\xi$ is the Poincaré dual class of a point on~$C$ and $\eta$ the mixed class defined above.
\end{proposition}

\subsection{Building a holomorphic line bundle}
\label{Ssec:HolomL}

We define a complex hermitian line bundle $L \to C$ by $L = \C \times_{U(1)} P$. For each point $w \in W$, the hermitian structure on~$L$ (inherited from $\C$), the complex structure on~$C$, and the connection $\nabla_t$ automatically endow~$L$ with the structure of a holomophic line bundle via the operator $\dbar = \nabla_w^{0,1}$. This holomorphic structure depends on~$w$. The quotient $T = W / H^1(C,\Z)$ is canonically identified with the $d$-th component of the Picard group of~$C$ parametrizing holomorphic structures on degree~$d$ line bundles:
$T = \rPic^d(C)$.

Recall that $H^0(C,\omega)$ is the space of holomorphic 1-forms on~$C$. It is a complex vector space of dimension~$g$. Its dual vector space $[H^0(C,\omega)]^*$ can be realized as the space of antiholomorphic 1-forms. 
We identify the real vector space $H^1(C,\R)$, which is the underlying vector space of~$W$, with $[H^0(C,\omega)]^*$ by assigning to an antiholomorphic 1-form $\alpha$ its real part:
$$
\begin{array}{rcl}
[H^0(C,\omega)]^* & \to & H^1(C,\R) \\
\mbox{antiholomorphic }\alpha& \mapsto & \Re(\alpha).
\end{array}
$$
Thus $W$ becomes a complex affine space. As before, we introduce the line bundle $L = P \times_{U(1)} \C \to W$ and define $\dbar = \nabla^{0,1}$ to endow it with a holomorphic structure. We now need to check that $\dbar^2=0$ or, equivalently, that the $(0,2)$-part of the curvature of $\nabla$ vanishes. This condition explains the choice of using antiholomorphic 1-forms rather than the holomorphic ones to endow~$W$ with a complex structure. Indeed, if $\alpha$ is an antiholomorphic 1-form, then it is itself the $(0,1)$-part of~$\nabla$. Then the $(0,2)$-part of the curvature is 
$$
\frac12\frac{\partial \alpha}{\partial \bar{\alpha}} = 0.
$$
Thus $\dbar$ endows the line bundle $L \to C \times W$ with a holomorphic structure. Taking the quotient by $H^1(C,\Z)$ we obtain a holomorphic line bundle over $C \times \rPic^d(C)$.

\begin{remark}
Since $\rPic^d(C)$ parametrizes holomorphic degree~$d$ line bundles over~$C$ it is natural to expect the existence of a universal line bundle over $C \times \rPic^d(C)$. This line bundle, however, is only defined up to a tensor product with the pull-back of any holomorphic line bundle over $\rPic^d(C)$. Under a tensor product like that the line bundle over each fiber $C \times w$ gets tensored by a line; thus its isomorphism class remains unchanged. 
The construction of the previous section reduces this ambiguity: the line bundle~$L$ is defined up to the pull-back of a {\em degree~0} holomorphic line bundle over $\rPic^d(C)$. Tensoring by a degree~0 line bundle corresponds to changing the gluing phases in the circle bundle. In particular the first Chern class of~$L$ (that we computed in the previous section) is well-defined in cohomology or in the Neron-Severi group, but not in Chow. 
\end{remark}

\subsection{One particle ground states}
\label{Ssec:1particle}

First consider one spinless particle on~$C$ in the presence of a magnetic field~$B$ and no electric field. The wave function for this particle is a square-integrable section of the line bundle $L \to C$. The Hamiltonian reads
$$
H \psi = \left(-\frac12\Delta_\nabla + a_1 B + a_2 R\right) \psi,
$$
where $\Delta = \nabla^* \nabla$ is the Laplacian of the connection, $R$ is the curvature of the Riemannian metric on the surface, and $a_1$, $a_2$ are some real constants that depend on the material. If charged particles have spin, the above Hamiltonian can still be a good approximation to determine the ground states, provided the magnetic field is strong enough to force the spin of the particle to be aligned with it.

The Weitzenb\"ock identity~\cite[Ch.~7, Cor.~(1.3)]{Demailly},  relates\footnote{In Demailly's notation, $\Delta''$ is our $\Delta_\dbar$, while $\Delta$ is our $\Delta_\nabla$. We have $\Delta = \Delta'+ \Delta''$ by Eq.~(1.15), so that
$$
\Delta = 2 \Delta'' - [i \Theta, \Lambda],
$$
where $i \Theta$ is the operator of multiplication by the curvature (in our case $B \omega$) and $\Lambda$ is the adjoint operator of multiplication by $\omega$ (in our case, simply dividing a 2-form by $\omega$).}
the Laplacian $\Delta_\nabla$ with the Laplacian $\Delta_\dbar = \dbar^* \dbar$:
$$
\Delta_\nabla = 2 \Delta_{\dbar} + B.
$$
Thus we can rewrite the Hamiltonian as
$$
H = - \Delta_\dbar + (a_1-1/2)B + a_2 R.
$$
Now we {\em require} that $(a_1-1/2)b + a_2 R$ be a constant function on~$C$. This can be reasonably achieved in different ways: one may assume that $a_2$ is negligibly small and $B$ is constant, or that both $B$ and $R$ are constant, or, more generally, one can adjust $B$ to compensate for the variation of~$R$. Whatever the choice, one obtains the Hamiltonian
$$
H = - \Delta_\dbar + a,
$$
where $a$ is a constant. The ground states of this Hamiltonian, i.e., the eigenfunctions $\psi$ with lowest eigenvalue, are the holomorphic sections of~$L$, because for a holomorphic section $\psi$ we have $\Delta_\dbar \psi = \dbar^* \dbar \psi = 0$. Note that the lowest level eigenspace only depends on the complex structure of~$C$.

\begin{remark}
The observation that the lowest eigenspace of~$H$ is highly degenerate is due to Landau. In some cases, up to a certain energy, the higher eigenstates are either also highly degenerate, or are organized in bands of eigenstates with close eigenvalues separated by large gaps proportional to~$d$. These bands, when present, are called {\em Landau levels}. The ground states of~$H$ form the lowest Landau level. If the charged particles are noninteracting fermions, at low temperature they will fill the Landau levels one by one.
\end{remark}

Thus the space of ground states of the one-particle Hamiltonian is $H^0(C,L)$, the space of holomorphic sections of~$L$. Its dimension $h^0(C,L)$ is the de\-ge\-ne\-racy of the ground state. It can be found via the Riemann-Roch formula:
$$
h^0(C,L) - h^1(C,L) = d + 1 -g.
$$
The annoying extra term $h^1(C,L)$ vanishes whenever $d > 2g-2$ (this is the simplest instance of the so-called Kodaira vanishing). Thus, if the total magnetic flux through the surface exceeds $2g-2$, then the degeneracy of the ground state equals $d+1-g$.

The purpose of this paper is to extend this method to $N$ interacting particles in the presence of the electric field. This brings on several changes. First, instead of a line bundle over~$C$ we will have a line bundle over the $N$th symmetric power of~$C$. The holomorphic sections of this line bundle are called {\em Laughlin states}. Second, instead of one line bundle~$L$ we consider the universal family of line bundles over $\rPic^d(C)$ constructed in the previous section. The Laughlin states form a vector bundle over $\rPic^d(C)$. The first Chern class of this vector bundle is closely related to the Hall conductance and can be computed by the Grothendieck-Riemann-Roch formula.

\begin{remark}
Laughlin's idea~\cite{Lau} to model quantum states of interacting charged particles with holomorphic wave-functions, or holomorphic sections of a line bundle, is actually more general: other models use sections of carefully chosen sheaves on the symmetric power of~$C$ (the Jain states, the Moore-Read states, and so on). We will not study them here, though our method potentially applies to some of them.
\end{remark}

\subsection{Laughlin states}
\label{Ssec:LaughlinStates}

Consider the $N$th power $C^N$ of the surface. On $C^N$, we define the line bundle $\L = \pi^*_1 L \otimes \pi^*_2 L \otimes \dots \otimes \pi^*_N L$. Here $\pi_i : C^N \to C$ is the projection to the $i$th factor.

Wave functions of a system of~$N$ indistinguishable particles on~$C$ are square integrable sections of~$\L$ over~$C^N$ that are either symmetric (for bosons) or antisymmetric (for fermions) with respect to the action of the symmetric group $S_N$ permuting the factors of $C^N$. A somewhat less precise way of saying the same thing is that the wave function is an (anti)symmetric function
$$
\psi(z_1, \dots, z_N)
$$
in some local coordinate $z$ on~$C$, such that if we fix all coordinates $z_j$ to some fixed values $z^{(0)}_j$ and only allow one coordinate $z_i$ to vary, the one-variable wave function 
$$
\psi \left(z^{(0)}_1, \dots, z^{(0)}_{i-1},\; z_i,\; z^{(0)}_{i+1}, \dots, z^{(0)}_N\right)
$$
is actually a section of~$L$ for any~$i$.

By definition, {\em Laughlin states} are (anti)symmetric holomorphic sections of $\L$ over $C^N$ characterized by one extra condition: they vanish to the order~$b$ at the diagonals $z_i=z_j$. Here $b \geq 1$ is an integral parameter; $b$ is even for bosons and odd for fermions.
Denote by $\Delta = \bigcup_{i < j} \{z_i = z_j \}$ the union of diagonals in $C^N$. It is a reducible divisor (that is, a linear combination of hypersurfaces). In more mathematical terms, a Laughlin state is an $S_N$-invariant holomorphic section of $\L(-b \Delta)$ for an appropriate lifting of the $S_N$ action to the line bundle.

As we will see in more detail in Proposition~\ref{Prop:S^NC}, the symmetric power of the surface
$$
S^NC = C^N / S_N
$$
is itself a smooth algebraic variety. Moreover, it carries a line bundle whose pull-back under the projection $C^N \to S^NC$ is $\L(-b\Delta)$. We will denote this line bundle by $\L_b$. Then a Laughlin state is simply a holomorphic section of $\L_b$ over~$S^NC$. This is the definition we will use. 

In the original paper~\cite{Lau}, Laughlin worked with $N$~particles on the plane with a constant magnetic field $B$. He wrote the wave function in the form
\begin{equation}\label{Laugh}
\Psi(z_1,...,z_N)=P(z_1,...,z_N) \prod_{i<j}(z_i-z_j)^b\cdot e^{-\frac B2\sum_{i=1}^N|z_i|^2 },
\end{equation}
where $z_1, \dots z_N$ are the coordinates of the particles, and $P(z_1,...,z_N)$ is a symmetric polynomial. 

Let's explain how this expression relates to our definition of Laughlin states. The first factor,
$$
P(z_1,...,z_N) \prod_{i<j}(z_i-z_j)^b,
$$
is a holomorphic section of the trivial line bundle over~$\C$ vanishing to the order~$b$ along the diagonals. This is precisely what we call a Laughlin state, except that $\C$ is not a compact Riemann surface. Because of the noncompactness, the space of holomorphic sections is infinite-dimensional, and the polynomial~$P$ could have been replaced by any symmetric holomorphic function, but Laughlin chooses polynomials, because they are eigenstates of the angular momentum operator. The second factor,
$$
e^{-\frac B2\sum_{i=1}^N|z_i|^2},
$$
is the square root of the hermitian metric~$h$ on the trivial line bundle over $\C$ whose curvature equals~$B$. It is present to ensure that $|\Psi|^2$ is the probability density for the particles distribution on the plane. In our more general setting we do not multiply by this factor, since it involves choosing local coordinates. Instead we simply write $|\psi|^2_h$ for the probability density.

The logic behind the definition of the Laughlin states is two-fold. A realistic hamiltonian for $N$ particles in a magnetic field with two-body interactions will have a generic form
\begin{equation}
H=-\sum_{i=1}^N\-\Delta_{\dbar,i} +\sum_{i<j}^NV(z_i,z_j)+\sum_{i=1}^N W(z_i),
\end{equation} 
where $\Delta_{\dbar,i}$ is the magnetic laplacian for the $i$th particle, $V: S^2C \to \R$ is the potential of the interaction and $W: C \to \R$ the potential related to impurities of the material.
Requiring the holomorphicity of the wave-function nullifies the first term (called the kinetic term). 
The order-$b$ vanishing condition minimizes the repulsive interactions, such as Coulomb forces. 
Moreover, the impurity potential $W$, if localised around some points $w_1,...,w_m$ on the surface, can be modelled by replacing the line bundle $L$ with $L(w_1+...+w_m)$.

%There exists a somewhat artificial choice of interaction, for which the Laughlin states are exact ground states, namely
%\begin{equation}
%V(z_i,z_j)=\sum_{\alpha=1}^b\bar\partial^\alpha_{z_i}\delta_{z_j}(z_i)%\partial^\alpha_{z_i},
%\end{equation}
%where $\delta_{z_j}$ is Dirac delta at $z_j$ and $\partial_{z}$ is a holomorphic derivative in local coordinates around $z_j$.

The main point here is that, as we will mention in the next section, if the Hamiltonian possesses a vector bundle of ground states separated from all other eigenstates by an energy gap, the Hall conductance can be recovered from the first Chern class and the rank of this vector bundle with no other information about the Hamiltonian. Deciding whether a particular Hamiltonian's ground states are separated by an energy gap is well-known hard question that we do not touch in this paper. 

The vector bundle of Laughlin states is assumed to be close enough to the vector bundle of ground states of the actual Hamiltonian to share the same rank and characteristic classes. Thus the results that we obtain for Laughlin states in this paper are actually valid for a wide range of Hamiltonians. On the other hand, the advantage of the Laughlin vector bundle is that it only depends on the holomorphic structure of~$C$ and~$L$ and therefore can be studied by methods of algebraic geometry.

\subsection{The vector bundle of Laughlin states}
\label{Ssec:Laugh}

In Section~\ref{Ssec:HolomL} we constructed a universal holomorphic line bundle $L$ over $C \times \rPic^d(C)$ up to a tensor product with a degree~0 line bundle over $\rPic^d(C)$. This construction upgrades in a straightforward way to a line bundle $\L_b$ over $S^NC \times \rPic^d(C)$. Denote by 
$$
\pi: S^NC \times \rPic^d(C) \to \rPic^d(C)
$$
the second projection. There is a vector bundle over $\rPic^d(C)$, denoted by $V=R^0\pi_*(\L_b)$ whose fiber $V_w$ over a point $w \in \rPic^d(C)$ is the corresponding space of Laughlin states. See Section~\ref{Ssec:Kodaira} for more mathematical details on the construction. We call it the {\em vector bundle of Laughlin states} or just the {\em Laughlin vector bundle}.

\subsection{The electric current in quantum mechanics}

\paragraph{One particle on the line.}
If $\psi(x)$ is the wave function of a charged particle on the line, the electric current operator is the same as the momentum operator:
$-i \d/\d x$ (where, of course, we have set $\hbar$, the mass, and the charge of the particle to~1). Note that, by the uncertainty principle, the current is not localized; for instance, it does not make sense to measure the electric current through a given point on the line, but only the global current along the line. Even an attempt to measure the current within some interval will fail. Indeed, it is easy to see that the operator
$$
\psi \to -i \rho(x) \d \psi / \d x
$$
is self-adjoint (hermitian) iff $\rho$ is a real constant function.

\paragraph{One particle on a surface.}
For a charged particle on a surface~$C$, any (real) tangent vector field~$T$ on~$C$ defines an operator
$$
I_T :\psi \to -i \nabla_T\psi,
$$
where $\nabla_T\psi$ is the covariant derivative of~$\psi$ along~$T$. Denote by $\beta = i_T \omega$, where $\omega$ is the symplectic form on~$C$. Then the operator $I_T$ is self-adjoint iff $\beta$ is closed. 

Indeed,
let $\psi_1$ and $\psi_2$ be two sections of~$L$. Denote by $(\cdot, \cdot)$ the pairing by the hermitian metric. We would like to check whether
$$
\left< I_T \psi_1, \psi_2 \right> - \overline{\left< I_T \psi_2, \psi_1 \right>}
=
-i \int_C \left[(\nabla_T \psi_1, \psi_2) + (\psi_1, \nabla_T \psi_2)\right] \omega
$$
vanishes for all $\psi_1, \psi_2$. 

We have
$$
i_T d(\psi_1, \psi_2) = (\nabla_T \psi_1, \psi_2) + (\psi_1, \nabla_T \psi_2),
$$
because the connection is compatible with the hermitian metric. On the other hand, $d(\psi_1, \psi_2) \wedge \omega = 0$ for dimension reasons. Hence
$$
0 = i_T [d(\psi_1, \psi_2) \wedge \omega]
= i_T d(\psi_1, \psi_2) \cdot \omega - d(\psi_1, \psi_2)  \wedge i_T \omega.
$$
We conclude that 
$$
\left[(\nabla_T \psi_1, \psi_2) + (\psi_1, \nabla_T \psi_2)\right] \omega = d(\psi_1, \psi_2)  \wedge i_T \omega = d(\psi_1, \psi_2) \wedge \beta.
$$
It follows that
$$
\left< I_T \psi_1, \psi_2 \right> - \overline{\left< I_T \psi_2, \psi_1 \right>}
= -i \int_C d(\psi_1, \psi_2)  \wedge \beta
= i \int_C (\psi_1, \psi_2) \wedge d\beta.
$$
This vanishes for any $\psi_1, \psi_2$ iff $d\beta=0$.

We have established that there is a current operator $I_T$ for any vector field~$T$ such that $\beta = i_T\omega$ is closed. If $\beta$ happens to also be exact, $\beta = df$, the corresponding operator measures the so called {\em circulating part} of the current. The case of interest for the quantum Hall effect, however, is when $\beta$ is a harmonic 1-form. Such 1-forms are in a one-to-one correspondence with 1-cohomology classes of~$C$ and measure the current flowing through the closed 1-cycles on the surface.

\begin{remark}
Let $T$ be a tangent vector field to~$C$. Then $\beta = i_T\omega$ is a harmonic 1-form iff $\alpha = i_Tg$ is. Indeed, we have $\alpha = \beta^*$, where $*$ is the Hodge star operator which, for 1-forms on Riemann surfaces, takes the form $(udx+vdy)^* = udy - vdx$ in any holomorphic local coordinate $z = x+iy$. The condition $\Delta \beta = 0$ is equivalent to $d\beta = d\beta* = 0$.

We use this remark in the next proposition.
\end{remark}

\begin{proposition}
Let $L \to C$ be a hermitian bundle with a hermitian connection~$\nabla$ over a surface~$C$ endowed with a riemannian metric~$g$. Let $T$ be a tangent vector field such that $\alpha = i_Tg$ and $\beta = i_T\omega$ are harmonic. Denote by $\nabla(s) = \nabla + s \alpha$, by $\Delta(s)$ the corresponding Laplacian and by $H(s) = -\frac12 \Delta_s + a_1B+a_2R$ the corresponding Hamiltonian. Then we have
$$
\nabla_T = \left. \frac{d H(s)}{ds}\right|_{s=0},
$$
so that the operator of measuring the electric current along~$T$ is equal to 
$$
-i \left. \frac{d H(s)}{ds}\right|_{s=0}.
$$
\end{proposition}

\paragraph{Proof.} The last two terms of the Hamiltonian do not depend on~$s$, so we only need to differentiate the Laplacian. We use the formula for the covariant Laplacian $\Delta = \nabla^* \nabla$ in local coordinates in terms of the connection $\nabla$ and the metric $g_{kj}$ given e.g. in  of~\cite[Ex.~10.1.32]{Nicolaescu}:
$$
\Delta  
= - \sum_{k,j} \left[ g^{kj} \nabla_k \nabla_j \; +\;  \frac{1}{\sqrt{|g|}}\, \partial_{x^k} ( \sqrt{|g|} g^{kj}) \cdot \nabla_j \right].
$$

Consider a holomorphic local coordinate $z$ on the Riemann surface~$C$ and write $z= x+iy$. Then in coordinates $(x,y)$, the metric has the form
$$
g = \left( 
\begin{array}{cc}
c(x,y) & 0 \\
0 & c(x,y)
\end{array}
\right),
\qquad
g^{-1} = \left( 
\begin{array}{cc}
1/c(x,y) & 0 \\
0 & 1/c(x,y)
\end{array}
\right),
$$
where $c(x,y)$ is a positive function. Note that $\sqrt{|g|} = c(x,y)$, so that the matrix $\sqrt{|g|} g^{-1}$ is just the identity matrix indepedently of $x,y$. Hence we have
$$
\partial_x (\sqrt{|g|} g^{-1}) = \partial_y (\sqrt{|g|} g^{-1})= 0.
$$
This kills the second term in the formula for the Laplacian drastically simplifying the expression: in coordinates $(x,y)$ as above we simply have
$$
\Delta = - \frac1{c(x,y)} [ \nabla_x^2 + \nabla_y^2].
$$
In coordinates we have $\nabla = d + \alpha^0$ and hence
\begin{align*}
\nabla_x(s) &= \partial_x + \alpha^0_x + s \alpha_x,\\
\nabla_y(s) &= \partial_y + \alpha^0_y + s \alpha_y,\\
\end{align*}
$$
\left. \frac{\partial \Delta(s)}{ds} \right|_{s=0}
= -\frac1{c(x,y)} \; \left[ 2 \alpha^0_x\alpha_x + 2 \alpha^0_y \alpha_y + \frac{\partial \alpha_x}{\partial x} + \frac{\partial \alpha_y}{\partial y} + 2\alpha_x \partial_x + 2\alpha_y \partial_y \right].
$$
and hence
$$
\left. \frac{\partial \Delta(s)}{ds} \right|_{s=0} = 
- \frac1{c(x,y)} \left[ 2 \alpha_x \nabla_x + 2 \alpha_y \nabla_y + \frac{\partial \alpha_x}{\partial x} + \frac{\partial \alpha_y}{\partial y} \right].
$$

Now, since $\alpha$ is a harmonic 1-form, we have
$$
 \frac{\partial \alpha_x}{\partial x} + \frac{\partial \alpha_y}{\partial y} = 0.
$$
(Proof: the form $\beta = \alpha^* = \alpha_y dx - \alpha_x dy$ is closed.) Moreover, since $\alpha$ was defined as $i_Tg$, we have
$$
\alpha_x = c(x,y) T_x, \quad \alpha_y = c(x,y) T_y.
$$
We finally get
$$
\left. \frac{\partial \Delta(s)}{ds} \right|_{s=0} = 
- 2 \left[ T_x \nabla_x + T_y \nabla_y \right] = -2 \nabla_T
$$
and
$$
\left. \frac{\partial H(s)}{ds} \right|_{s=0}  = \nabla_T.
$$
\qed

\paragraph{Many particles on a surface.}
As before, we take a tangent vector field $T$ on the surface~$C$ such that the 1-form $\alpha = i_Tg$ is harmonic. We denote by $\nabla(s)$ the connection $\nabla+s\alpha$ and by $H(s)$ the corresponding Hamiltonian.

For $N$ particles, the current operator $J_T$ is the sum of one-particle currents. As for the Hamiltonian, it is the sum of Hamiltonians of individual particles plus an interaction term that does not depend on~$s$. Thus the equality 
$$
J_T= -i \left. \frac{d H(s)}{ds}\right|_{s=0}
$$
is still true and can more generally be taken as a definition of the current operator.

\subsection{The Hall conductance}

The most well-known fact about the integer quantum Hall effect is that the Hall conductance in units $e^2/h$ is equal to the first Chern class of the line bundle of ground states. For the fractional quantum Hall effect the ground state is degenerate, so we have a vector bundle instead of a line bundle, and the first Chern class is replaced with the slope $=$ the first Chern class divided by the rank. 

There are different ways to make a connection between the Hall conductance and the slope of the vector bundle of ground states, but none is particularly straightforward, because the connection only appears when one takes an average over a long time. We chose not not review the rather lengthy proofs; instead we refer to the strongest result that belongs to Avron and Seiler~\cite{AS}. Using the definition of the electric current as above, they proved in a mathematically rigorous way\footnote{Stricty speaking, they consider the genus~1 case, but we believe that the argument extends with very little change to any genus.} that for a constant electric field $\varepsilon E$, the expectation value of the Hall conductance averaged over the whole Picard group tends to the slope of this vector bundle faster than any power $\varepsilon$. This fast convergence rate probably explains why the quantum Hall experiments give such precise results.

\subsection{Prior state of the art in physics}
The integer QHE (when the Hall conductance equals $\nu \frac{e^2}{h}$ with an integer~$\nu$) admits a well-developed mathematical explanation in terms of non-interacting fermions, making use of various techniques, from spectral indices of compact operators, to non-commutative geometry to Chern classes on moduli space of flat connections, see, e.g.,~\cite{Av} for a review.

The fractional QHE (the case of a fractional $\nu$) is much harder to explain analytically as it corresponds to strongly-interacting electrons. Laughlin's approach, that consists in guessing a good approximation of the many-particles wave functions, explains the plateaux with simple fractions $\nu=1/b$, $b\in\mathbb N$. The idea to geometrize the problem goes back again to Laughlin, who considered the integer QHE on a cylinder. Later on, Haldane \cite{H1983}, Haldane-Rezayi \cite{HR} and Wen-Niu \cite{WN} developed the theory of Laughlin states in the fractional QHE, on the sphere, torus and compact surfaces of higher genus, respectively. 
Avron, Seiler and collaborators have developed the adiabatic transport of QHE states on Riemann surfaces \cite{AS,ASY,ASZ,KS}. As a result, the subject gave rise to a new field of topological phases of matter, with potential applications to quantum computers, novel materials, etc. 

There is a number of interesting physics conjectures with regard to Laughlin states and more general FQHE states on Riemann surfaces.  In this paper we prove some of these conjectures using the methods of algebraic geometry.

As follows from the discussion in Section~\ref{Ssec:Laugh}, the Laughlin vector bundle $V=V_{N,d,b,g}$ is indexed by the following positive integers --  the number of particles $N$, the degree of the magnetic line bundle $d$ (flux of the magnetic field through the surface), the  order of vanishing $b$ and the genus of the surface $C$. The fractional number $1/b$ is called the filling fraction.

\paragraph{A. The ``Topological shift" conjecture} (Wen-Zee \cite{WZ}, see also Fröhlich-Studer \cite{FS}).
For an incompressible quantum Hall state on a compact closed Riemann surface, there exists $\nu\in\mathbb Q$, such that $N,d,\nu$ and $g$ are related as
\begin{equation}
\nu^{-1}N=d+S,
\end{equation}
where the quantum number $S$ is called the shift and it is ``topological" in the sense that it depends only on $g$. This relation is known as Wen-Zee  formula.

For the Laughlin states we have $\nu^{-1} = b$ and $S = b(1-g)$, and the conjecture reads: the space of Laughlin states has positive dimension for
$bN = d+b(1-g)$ and dimension~0 for $bN > d+b(1-g)$. When $bN = d+b(1-g)$, 
the corresponding Laughlin states are commonly called ``completely filled" with filling fraction $\nu=1/b$. 

\paragraph{B. The ``Topological degeneracy" conjecture} (Haldane-Rezayi \cite{HR} for $g=1$, Wen-Niu \cite{WN} for any $g$).
For completely filled states defined as above, the rank of the Laughlin vector bundle equals
\begin{equation}
{\rm rk\,}V_{N,d,b,g}=b^g.
\end{equation}
Notably, it is ``topological" in the sense of being independent of $N$ and $d$. 

This conjecture is a foundational result of the so-called topological states of matter. A construction of $b^g$ linearly independent Laughlin states on a genus~$g$ surface was carried out explicitly in terms of the Riemann theta functions and prime forms in~\cite{SK2019}. This paper contains the first rigorous proof that this family is indeed a basis.

\paragraph{C. Hall conductance for completely filled states.}

In the case $bN = d+b(1-g)$, the Hall conductance equals $1/b$ times the 
Theta class on $\rPic^d(C)$, or, in terms of differential forms,
\begin{equation}
\sigma_H=\frac1b\left[\sum_{i=1}^g \alpha_i\wedge\beta_{i}\right],
\end{equation}
where $(\alpha_i,\beta_i)$ is a symplectic basis of $H^1(\rPic^d(C),\R)$.

For $g=1$ this observation goes back to \cite{NTW,TW,AS} and others, see also \cite{Var} and \cite{BK1}. Its analog for higher genus Riemann surfaces was discussed in the case of Integer QHE by Avron-Seiler-Zograf \cite{ASZ}, see also \cite{KMMW}. 

We also note that the Hall conductance on a torus was computed by Hasting-Michalakis for the gapped local Hamiltonian without making use of the Laughlin states~\cite{HM}.

\paragraph{D. The ``Projective flatness'' conjecture.}
This conjecture states that the vector bundles of completely filled quantum Hall states over geometric parameter spaces are projectively flat for an appropriate connection. %{\color{red} and $L^2$ structure}. 

This conjecture is implicit in early papers on quantum Hall states~\cite{WN} and is often discussed in the case of the moduli space of the complex structures, see e.g.\cite{KVW,Read,KW}. In~\cite{KZ2021}, we proposed the projective flatness criterion as a test for the topological nature of a quantum state of matter. For $g=1$, flat metric on the torus and the constant magnetic field the projective flatness of Laughlin was proved in \cite{BK2}.

For the introduction for some of the topics discussed we refer to~\cite{Av,SK2015}. We note also other approaches due to Belissard et al \cite{BES}, Froehlich et al \cite{FST}, Bachmann et al \cite{Bach} and others.

In the next section we describe our results, in particular, we prove Conjectures A, B and~C for Laughlin states and make a step towards conjecture~D by showing that the Chern character of the Laughlin bundle is compatible with projective flatness.

\subsection{Main mathematical results}

The physical introduction has led us to the following mathematical setup. 
We fix a smooth algebraic curve~$C$ of genus~$g$ and a degree~$d$ line bundle $L$ over~$C$. The product $C^N \times \rPic^d(C)$ carries the line bundle
$$
\L(-b\Delta) = \pi_1^*L \otimes \dots \otimes \pi_N^*L 
\left(-b\sum_{i-j} \Delta_{ij}\right).
$$

\subsubsection{The line bundle descends to the symmetric product $S^NC$}

We lift the natural $S_N$ action on $C^N$ to the total space of the line bundle in the following way: if $b$ is even, the action lifts by permuting the factors $\pi_i^*L$; if $b$ is odd, the outcome is further multiplied by $\pm1$ for even/odd permutations. We are interested in the space of global $S_N$-invariant sections of $\L(-b\Delta)$ over~$C^N$. Denote by $S^NC$ the $N$th symmetric power of~$C$.

\begin{proposition} \label{Prop:S^NC}
There exists a line bundle $\L_b$ on $S^NC$ such that
(i)~its pull-back under the quotient map $q: C^N \to S^NC$ is equal to $\L(-b\Delta)$, and 
(ii) the space of its holomorphic sections over an open set $U \subset S^NC$ is equal to the space of $S_N$-invariant sections of $\L(-b\Delta)$ over $\pi^{-1}(U)$.
\end{proposition}

\subsubsection{The first Chern class of $\L_b$}

There is a natural Abel-Jacobi map $\sigma: S^NC \to \rPic^N(C)$ assigning to a collection of points the corresponding divisor up to rational equivalence. Thus $S^NC$ carries the cohomology class $\theta_N \in H^2(S^NC,\Z)$: the pull-back from $\rPic^N(C)$ of the cohomology class Poincaré dual to the $\Theta$-divisor. Further, it is proved in~\cite{Mattuck}, that for $N > 2g-1$ the map $\sigma$ is a projective bundle. In other words, there exists a vector bundle $E \to \rPic^N(C)$ such that $S^NC = \P(E)$. We denote the rank of~$E$ by
\begin{equation} \label{eq:r}
r = N-g+1.
\end{equation}
The geometry of $S^NC$ is more complicated for smaller values of~$N$. In this paper we assume that $N > 2g-1$, although some of our results might also hold for smaller values of~$N$.

The vector bundle~$E$ is defined up to a tensor product with a line bundle. After specifying one of the choices, one defines the first Chern class $\xi = c_1(\cO(1))$. 

\begin{proposition}
We have
$$
c_1(\L_b) = b \theta_N + p \xi,
$$
where 
\begin{equation} \label{Eq:p}
p= d- b(N+g-1).
\end{equation}
\end{proposition}

The integer~$p$, called the {\em number of quasi-holes} in physics, will play an important role in our formulas. In particular, it follows from the expression of the Chern class that the restriction of $\L_b$ to a fiber of~$\sigma$ is the line bundle $\cO(p)$ on a projective space $\CP^{r-1}$. It has no sections for $p<0$ and only constant sections for $p=0$. Its higher cohomology groups $H^i(\CP^{r-1}, \cO(p))$ vanish for $p > -r$.

\subsubsection{Extending the above to a universal family of line bundles}

Over $C \times \rPic^d(C)$ we consider a universal degree~$d$ line bundle~$L$ such that the restriction of $L$ to ${\rm point} \times \rPic^d(C)$ has zero first Chern class.
 
Both of the above propositions upgrade to this setting. The formula for the first Chern class of $\L_b$ acquires an extra term, see Proposition~\ref{Prop:c_1(L)}.

\subsubsection{The Laughlin bundle and its characteristic classes}

The two following theorems constitute the main result of the paper.

Let $\pi : S^NC \times \rPic^d(C) \to  \rPic^d(C)$ be the projection to the second factor. Denote by $\theta_d \in H^2(\rPic^d(C),\Z)$ the cohomology class Poincaré dual to the theta-divisor. 

\begin{theorem}[Vanishing]\label{thm:Kodaira}
Assume $N > 2g-1$ and $p > -r$. Then
the sheaves $R^i \pi_* \L_b$ on $\rPic^d(C)$ vanish for $i \geq 1$; the sheaf $R^0 \pi_* \L_b$ is locally free, that is, the sheaf of sections of a vector bundle.
\end{theorem}

We denote the vector bundle $R^0 \pi_* \L_b$ by $V$ or $V_{N,d,b,g}$ and call it the {\em Laughlin vector bundle}.

\begin{theorem}[Chern character]\label{thm:ch}
Still assuming $N > 2g-1$ and $p > -r$, the Chern characters of the Laughlin vector bundle over $\rPic^d(C)$ equal
$$
{\ch}_m(V)=\sum_{k=m}^g
{{g-m}\choose{k-m}}{{N-g+p}\choose{k-g+p}}b^{k-m} \frac{(-\theta_d)^m}{m!}.
$$
Here $\theta_d$ is the Theta class on $\rPic^d(C)$, and we use the convention that if $k-g+p <0$ in the second binomial coefficient, the corresponding term vanishes.
\end{theorem}

Let's consider some particular cases of Theorem~\ref{thm:ch}. 

\paragraph{The case $p<0$.} In this case, the line bundle $\L_b$ has no sections, so $V$ is the rank~0 vector bundle and $\ch(V) = 0$. As for the formula, all the terms vanish because $k - g+ p < 0$ for all~$k$. This proves the Wen-Zee formula (Conjecture A) for the Laughlin states with $S=b(1-g)$.

\paragraph{The case $p=0$.} In this case, the formula simplifies drastically. The only value of $k$ for which $k-g+p  = k-g$ is nonnegative is $k=g$. Thus we get
$$
\ch_m(V) = b^{g-m} \frac{(-\theta_d)^m}{m!}
$$
or, in other words,
$$
\ch(V) = b^g \exp(-\theta_d/b).
$$
In particular, the rank of~$V$ (the ``topological degeneracy") equals~$b^g$, which proves Conjecture~B, and the Hall conductance, given by the slope of $V=\frac{c_1(V)}{{\rm rk\,}V}$, equals $-\theta_d/b$, proving Conjecture~C.

Recall that (see~\cite[Ch II, \S 3]{Kob}), if $V$ is projectively flat, then $$\ch(V)={\rm rk\,} V\cdot \exp\left(\frac{c_1(V)}{{\rm rk\,} V}\right)$$.

Thus, we can show that the vector bundles with $p>0$ are not projectively flat, at least starting from some genus. At $p=0$, the characteristic classes of $V$ are consistent with being projectively flat (cf. Conjecture~D).

\paragraph{Genus~$0$.} In this case $\rPic^d(C)$ is a point, so that the Laughlin bundle is just a vector space. Its dimension equals
$$
{N+p \choose p}.
$$

\paragraph{Genus~$1$.} In this case $\rPic^d(C)$ is isomorphic~$C$, $\theta_d$ is Poincaré dual to the class of a point in~$C$,
$$
{\rm rk\,} V = {N+p-1 \choose p} \left(b +\frac{p}{N}\right),
$$
$$
c_1(V) = \ch_1(V) = -{N+p-1 \choose p} \theta_d.
$$
The slope of~$V$ equals
$$
\frac{c_1(V)}{{\rm rk\,} V} = -\frac{\theta_d}{b +p/N}.
$$

\subsection{Notation}

Here is some notation consistently used in the paper.

\begin{itemize}
\item $C$ a Riemann surface;
\item $g \geq 0$ its genus;
\item $L$ a line bundle over~$C$ or, by abuse of notation, a universal line bundle over $C \times \rPic^d(C)$;
\item $d$ its degree;
\item $N$ the number of particles;
\item $b \geq 1$ a positive integer;
\item $p = d - b(N+g-1)$ the number of quasi-holes;
\item $\rPic^k(C)$ the $k$th component of the Picard group of~$C$;
\item $\theta_k \in H^2(\rPic^k(C),\Z)$ the cohomology class in $\rPic^k(C)$ Poincaré dual to the Theta divisor;
\item $S^NC$ the $N$th symmetric power of~$C$;
\item $q: C^N \to S^NC$ the quotient map;
\item $\sigma: S^NC \to \rPic^N(C)$ the Abel-Jacobi map;
\item $\pi_i C^N \to C$ the $i$th projection;
\item $\pi : S^NC \times \rPic^d(C) \to \rPic^d(C)$ the second projection;
\item $\Delta$ the union of diagonals in $C^N$;
\item $\L = \pi_1^*L \otimes \dots \otimes \pi_N^*L$ a line bundle over $C^N$ or, by abuse of notation the corresponding universal line bundle over $C^N \times \rPic^d(C)$;
\item $\L_b$ the line bundle over $S^NC$ whose pull-back to $C^N$ is $\L(-b\Delta)$ or, by abuse of notation, the corresponding universal line bundle over $S^NC \times \rPic^d(C)$.
\end{itemize}

\section{Cohomology and Chow ring of $S^NC$ (after A.~Mattuck)}

In this section we summarize well-know facts on $\Theta$-divisors and A.~Mattuck's results from~\cite{Mattuck} on $S^NC$.

\subsection{$\Theta$-divisors}

We denote by $\Theta$ the canonical $\Theta$-divisor in $\rPic^{g-1}(C)$. The Chow classes of its translates differ from each other; thus the Chow group of $\rPic^{g-1}(C)$ is not finitely generated. However, the translations act trivially on the cohomology group of $\rPic^{g-1}(C)$. Thus the Poincaré dual cohomology class is the same for any translate of $\Theta$.

The Picard group components $\rPic^{g-1}(C)$ and $\rPic^k(C)$ for any $k \in \Z$ are naturally identified up to translation. Hence there is a well defined cohomology class $\theta_k \in H^2(\rPic^k(C),\Z)$ for  every $k \in \Z$. If we choose a symplectic basis $(\alpha_i, \beta_i)_{1 \leq i \leq g}$ of $H^1(\rPic^k(C),\Z)$, this class is equal to 
$$
\theta_k = \sum_{i=1}^g \alpha_i \wedge \beta_i.
$$

For an integer $b \geq 1$, consider the map of multiplication by~$b$ 
$$
\mathbf{b}: \rPic^k(C) \to \rPic^{bk}(C).
$$
We have 
\begin{equation} \label{Eq:pullbackb}
\mathbf{b}^*(\theta_{bk}) = b^2 \theta_k.
\end{equation}

Let $(\alpha_i,\beta_i)_{i=1}^g$ be a symplectic basis of $H^1(\rPic^k(C),\Z)$ and $(\alpha'_i,\beta'_i)_{i=1}^g$ a symplectic basis of $H^1(\rPic^{k'}(C),\Z)$.

\begin{notation} \label{Not:eta}
In the product $\rPic^k(C) \times \rPic^{k'}(C)$ denote 
by $\eta$ the ``mixed'' cohomology class
$$
\eta_{k,k'} = \sum_{i=1}^g (\alpha_i \wedge \beta'_i  +  \alpha'_i \wedge \beta_i).
$$
\end{notation}

It is easy to see that this class does not depend on the order of the factors or on the choice of the symplectic basis. We will suppress the indices $k$ and $k'$ in the notation if this causes no confusion.

\begin{remark}
We use the same notation as in Proposition~\ref{Prop:c1OneParticle}; indeed if we restrict $\eta$ to $\rPic^d(C) \times C$, where $C \subset \rPic^1(C)$ is the standard Abel-Jacobi embedding of the curve into its Jacobian, the class $\eta$ defined above restricts to the class~$\eta$ from Proposition~\ref{Prop:c1OneParticle}.
\end{remark}

Consider the summation map
$$
\mathbf{s}: \rPic^k(C) \times \rPic^{k'}(C) \to \rPic^{k+k'}(C).
$$
The pull-back of $\theta_{k+k'}$ under this map equals
\begin{equation} \label{Eq:pullbacks}
\mathbf{s}^* \theta_{k+k'} = \theta_k + \theta_{k'} + \eta.
\end{equation}
This formula is obtained simply by expanding 
$$
\sum_{i=1}^g (\alpha_i+\alpha_i') \wedge (\beta_i + \beta'_i).
$$
The pull-back of $\eta$ under the multiplication map 
$$
\beta: \rPic^k(C) \times \rPic^{k'}(C) \to \rPic^{bk}(C) \times \rPic^{b'k'}(C)
$$
equals
\begin{equation} \label{Eq:pullbacketa}
\beta^*\eta_{bk,b'k'} = bb' \eta_{k,k'}.
\end{equation}
Finally, we will prove an imporant equality in the Chow ring.
Recall that the $\Theta$-divisor in $\rPic^{g-1}(C)$ is the locus of line bundles with a nonzero section.

\begin{lemma} \label{Lem:degg}
Let $S$ be a degree~$g$ line bundle over a curve~$C$ and consider the map 
$$
\begin{array}{ccccc}
\phi&:&C &\to & \rPic^{g-1}(C)\\
&& z & \mapsto & S(-z).
\end{array}
$$
Then we have a non-canonical isomorphism of line bundles $\phi^*\cO(\Theta) \isom S$, in other words, $c_1(\phi^*\cO(\Theta)) = c_1(S)$ in the Chow ring of~$C$.
\end{lemma}

\paragraph{Proof.} First assume that $S$ is generic in the sense that it has a unique nonzero section $s$ up to multiplicative constant and this section has~$g$ simple zeros. Then the pull-back $\phi^{-1}(\Theta)$ is a degree~$g$ divisor on~$C$ composed precisely of the zeros of~$s$. Indeed, $S(-z)$ lies in the theta-divisor iff $S(-z)$ has a nonzero section, which is equivalent to $z$ being a zero of~$s$. We see that the divisor of $S$ and $\phi^{-1}(\Theta)$ coincide. Therefore $S \isom \phi^* \cO(\Theta)$. In other words, the map $\rPic^g(C) \to \rPic^g(C)$ that maps $S \mapsto \phi^* \cO(\Theta)$ is the identity map on the locus of generic divisors. By continuity, it is the identity map on the whole $\rPic^g(C)$, so $S = \phi^*\cO(\Theta)$ holds for any divisor~$S$. \qed

\subsection{The symmetric power $S^NC$}
\label{Ssec:S^NC}

Consider the natural map $\sigma: S^NC \to \rPic^N(C)$ assigning to an unordered collection of points the corresponding divisor up to rational equivalence. The fibers of this map are projective spaces: the linear systems corresponding to each divisor. For $N \geq 2g-1$ all these projective spaces have the same dimension~$N-g$.

Now we have to make a noncanonical choice, namely, fix a point $z_0 \in C$.

\begin{notation} \label{Not:xi}
Denote by~$\xi \subset S^NC$ the locus of collections of~$N$ points containing~$z_0$. In other words, $\{z_1, \dots, z_N \} \in \xi$ iff $z_0 \in \{z_1, \dots, z_N \}$. By abuse of notation, we will also denote by $\xi$ its Poincaré dual Chow or cohomology class in $A^1(X)$ or $H^2(S^NC,\Z)$.
\end{notation}

\begin{proposition}
The Chow class $\xi^k$ in $S^NC$ is represented by the cycle $\Xi_k$ of collections of $N$ points at least $k$ of which are equal to~$z_0$.
\end{proposition}

\paragraph{Proof.} Denote by $X_i \in C^N$ the locus of $N$-tuples $(z_1, \dots, z_N)$ such that $z_i = z_0$. Then $X_1 + \dots + X_N$ is the full preimage of $\xi$ in~$C^N$. Since the degree of the quotient map $q: C^N \to S^NC$ equals $N!$, we have
$$
\xi^k = \frac1{N!} q_* (X_1 + \dots + X_N)^k.
$$
Now, the normal vector bundle to $X_i$ is trivial, so that $X_i^2 = 0$. Therefore
$$
(X_1 + \dots + X_N)^k = k! \sum_{i_1< \dots < i_k} X_{i_1} \dots X_{i_k}.
$$
The right-hand side contains $\frac{N!}{k!\,(N-k)!}$ terms. In each term
the cycle $X_{i_1} \dots X_{i_k}$ maps onto $\Xi_k$ with degree $(N-k)!$. Putting everything together, we get
$$
\frac1{N!} \cdot k! \cdot \frac{N!}{k!\,(N-k)!} \cdot (N-k)! \cdot \Xi_k = \Xi_k.
$$
\qed

\begin{example}
If $C = \CP^1$, we have $S^NC = \CP^N$ and $\Xi_k = \CP^{N-k}$.
\end{example}

\begin{notation}\label{Not:W}
Further, for $0 \leq i \leq g$, denote by $W_i = \sigma_*\xi^{N-g+i} \subset \rPic^N(C)$ the locus of divisors in~$C$ that can be represented as 
$$
(N-g+i)z_0 + z_1 + \dots + z_{g-i},
$$
in other words, as $N-g+i$ times $z_0$ plus $g-i$ any other points. Once again, by abuse of notation, we also denote by $W_i$ the Poincaré dual cohomology class in $A^i(\rPic^N(C))$ or $H^{2i}(\rPic^N(C),\Z)$. We will also use the same notation for the pullback of $W_i$ to $S^NC$.
\end{notation}

We will be mostly interested in cohomology classes, since our goal is to apply the Grothendieck-Riemann-Roch formula, which only uses cohomology classes rather than Chow classes. However, for completeness we specify when an identity holds in the Chow ring.

Here are some equalities that only hold (or only make sense) in cohomology:
\begin{eqnarray} \label{Eq:W}
W_i &=& \frac{\theta_N^i} {i!} \, , \\
\nonumber
W_{g-1} &=& \mbox{Abel-Jacobi embedding of } C,\\
\nonumber
W_g & = & \mbox{point}.
\end{eqnarray}
In particular, the intersection index of the $\Theta$-divisor $W_1$ with the curve $W_{g-1}$ equals~$g$ by the following computation:
$$
W_1 \cdot W_{g-1} = \theta_N \cdot \frac{\theta_N^{g-1}}{(g-1)!} = g\,\, \frac{\theta_N^g}{g!} = g W_g = g [{\rm point}].
$$

\begin{theorem} [\cite{Mattuck}, Thm 3 and corollary]
Assume that $N \geq 2g-1$.
There exists a vector bundle $E \to \rPic^N(C)$ of rank $r = N-g+1$ such that 
\begin{itemize}
\item
$\P(E) = S^NC$, 
\item 
the first Chern class $c_1(\cO(1))$ on $\P(E)$ equals $\xi$ in the Chow group $A^1(S^NC)$.
\end{itemize}
\end{theorem}

\begin{corollary}
The total Segre class of~$E$ equals
$$
s(E) = \pi_* \sum_{k \geq 0} \xi^k = W_0 + W_1 + \dots + W_g
$$
in the Chow ring of~$\rPic^N(C)$ and
$$
s(E) = e^{\theta_N}
$$
in $H^*(\rPic^N(C),\Q)$. Hence the total Chern class of~$E$ equals
$$
c(E) = 1/s(E) = e^{-\theta_N}
$$
in $H^*(\rPic^N(C),\Q)$. In particular, we have
%$$
%A^*(S^NC) = A^*(\rPic^N(C)) [\xi] / \left(\xi^g - W_1 \xi^{g-1} + \dots +(-1)^g W_g\right)
%$$
%and
$$
H^*(S^NC, \Z) = H^*(\rPic^N(C),\Z) [\xi] / \left(\xi^r - \frac{\theta_N}{1!} \xi^{r-1} + \dots +(-1)^g \xi^{r-g} \frac{\theta_N^g}{g!}\right).
$$
\end{corollary}

\paragraph{Proof.} We have used two well-known facts on vector bundles $E \to B$ over an algebraic variety~$B$. Let $\xi = c_1(\cO(1)) \in H^2(B,\Z)$. Then the total Chern class of~$E$ and the total Segre class of~$E$ are inverse of each other:
$$
s(E): = \sigma_*(1 + \xi + \xi^2 + \dots) = 1/c(E).
$$
And the cohomology ring of $\P(E)$ is equal to the quotient of the polynomial ring $H^*(B,\Z)[\xi]$ by the Chern polynomial applied to $\xi$:
$$
\xi^r + c_1(E) \xi^{r-1} + \dots + c_r(E).
$$
\qed

%\begin{remark}
%The cohomology classes $W_i$ and $\xi$ do not depend on the choice of the point~$z_0$. The Chow classes, on the other hand, do depend on $z_0$, so that 
%$$
%\xi^g - W_1 \xi^{g-1} + \dots +(-1)^g W_g = 0
%$$
%is just one identity in cohomology, but a family of identities in the Chow ring.
%\end{remark}

\section{The first Chern class of $\L_b$} \label{Sec:Lb}

The existence of the line bundle $\L_b$ was announced in Proposition~\ref{Prop:S^NC}.

\paragraph{Proof of Proposition~\ref{Prop:S^NC}.} Consider a point $(z_1, \dots, z_N) \in C^N$ and let $I_1 \sqcup I_2 \sqcup \dots \sqcup I_k = \{ 1, \dots, N \}$ be the set partition into groups of coinciding points. For each part $I$ of the set partition, denote by $w_I$ a local coordinate on~$C$ in the neighborhood of the corresponding coinciding points. Let 
$$
V_I = \prod_{\substack{i,j \in I, \\ i<j} }  (w_I(z_i) - w_I(z_j))
$$
be the Vandermonde determinant. Finally, let $e_I$ be the collection of the elementary symmetric functions of $w_I(z_i)$ for $i \in I$. Then the $e_I$'s taken together form a family of local coordinates on $S^NC$ at the point 
$\{z_1, \dots, z_N\}$ and the sheaf of $S_N$-invariant sections of $\L(-b\Delta)$ in these coordinates is generated by the section $\prod_{j=1}^k V_{I_j}^b \times \alpha$, where $\alpha$ is any nonvanishing symmetric local section of~$\L$. Thus it is locally free of rank~1, and its generator pulls back to a generator of the sheaf of sections of the line bundle $\L(-b\Delta)$ on~$C^N$. This proves both statements of the proposition. \qed

Recall the notation for the integer~$p$ (Eq.~\eqref{Eq:p}) and for the 2-cohomology classes
\begin{itemize}
\item
$\xi$, see Notation~\ref{Not:xi};
\item $\theta_N$, Poincaré dual to the $\Theta$-divisor in $\rPic^N(C)$;
\item $\eta$, see Notation~\ref{Not:eta}.
\end{itemize}

\begin{proposition} \label{Prop:c_1(L)}
 The first Chern class of the line bundle $\L_b$ over $S^NC \times \rPic^d(C)$ equals 
$$
 b\theta_N - \eta + p \xi
 $$ 
 in $H^2(S^NC \times \rPic^d(C),\Z)$.
\end{proposition}

\begin{remark}
The cohomology classes $\theta_N$, $\eta$, and $\xi$ are defined on $\rPic^N(C)$, $\rPic^N(C) \times \rPic^d(C)$, and $S^NC$ respectively. In all three cases we take their pullbacks to $S^NC \times \rPic^d(C)$.
\end{remark}

\paragraph{Proof of the Proposition.} It is enough to check the statement on $C^N \times \rPic^d(C)$, because the map $H^2(S^NC,\Z) \to H^2(C^N,\Z)$ is an injection.

\paragraph{Basic case: $p=0$, $b=1$.} In this case our line bundle $L \to C$ has degree $d = N+g-1$. Consider the map
$$
\begin{array}{rccc}
s : & C^N \times \rPic^d(C) & \to & \rPic^{g-1}(C)\\
& (z_1, \dots, z_N, \quad L) & \mapsto & L(-\sum z_i).
\end{array}
$$

We claim that the line bundles 
$$
\L(-\Delta)
\qquad \text{and}
\qquad
s^*\cO(\Theta)
$$
over $C^N \times \rPic^d(C)$ are isomorphic up to a tensor product with a line bundle over $\rPic^d(C)$. To show that, fix $L \in \rPic^d(C)$ and the values of $z_1, \dots, z_{N-1}$ and let $z_N$ vary over~$C$. The restriction of the line bundle 
$$
\L(-\Delta)
$$ 
to this curve~$C$ is isomorphic to the degree~$g$ line bundle
$$
S = L \left(-\sum_{i=1}^{N-1} z_i \right).
$$
According to Lemma~\ref{Lem:degg} the restriction of the two line bundles to~$C$ are isomorphic, and this is true for any fiber of the projection 
$$
C^N \times \rPic^d(C) \to C^{N-1} \times \rPic^d(C). 
$$
It follows that the two line bundles are isomorphic up to a twist by a line bundle over~$C^{N-1} \times \rPic^d(C)$. To determine this twist we now fix $z_1, \dots, z_{N-2}$ and $z_N$ and let $z_{N-1}$ vary over~$C$. By the same token as above, the restriction of both line bundles to this curve are isomorphic. We conclude that the line bundles can only differ by a line bundle over $C^{N-2} \times \rPic^d(C)$. Performing this comparison $N-2$ more time we obtain that the two line bundles are isomorphic up to a twist by a line bundle over $\rPic^d(C)$. 

Now, the map $s$ is a composition of several maps that we have already studied: $\pi : S^N \to \rPic^N(C)$, the multiplication by $-1$ from $\rPic^N(C)$ to $\rPic^{-N}(C)$, and the summation map from $\rPic^{-N}(C) \times \rPic^d(C)$ to $\rPic^{g-1}(C)$. The pull-back of the cohomology class $\theta_{g-1}$ under the summation map equals $\theta_{-N} + \eta_{-N,d} + \theta_d$ by Equation~\eqref{Eq:pullbacks}. The pull-back of this under the multiplication by $-1$ equals $\theta_N - \eta + \theta_d$ by Equations~\eqref{Eq:pullbackb} and~\eqref{Eq:pullbacketa}. So the first Chern class of $\L(-\Delta)$ equals $\theta_N - \eta + \theta_d$ up to a pull-back from $\rPic^d(C)$. 

To lift this ambiguity, choose $N$ fixed pairwise distinct points $z_1, \dots, z_N \in C$ and consider the restriction of $\L(-\Delta)$ to $(z_1, \dots, z_N) \times \rPic^d(C)$. This restriction is the tensor product of $L|_{z_1} \otimes L|_{z_2} \otimes \dots \otimes L|_{z_N}$. As explained in Section~\ref{Ssec:Factorc1}, the first Chern class of $L$ restricted to $\mbox{point} \times \rPic^d(C)$ vanishes. Thus the restriction of $\L(-b\Delta)$ to $(z_1, \dots, z_N) \times \rPic^d(C)$ also vanishes. We conclude that 
$$
c_1(\L(-\Delta)) = \theta_N - \eta
$$
for $b=1$, $p=0$.

\paragraph{The general case.} Choose a point $z_0 \in C$ as in Section~\ref{Ssec:S^NC}. Our choice of universal line bundle $L$ over $C \times \rPic^d(C)$ determines a universal line bundle $L(-p z_0)$ over $C \times \rPic^{d-p}(C)$, where the shift by the divisor $-pz_0$ identifies $\rPic^d(C)$ and $\rPic^{d-p}(C)$. Now, by definition of~$p$, we have
$d-p = b (N+g-1)$, that is, the degree of $L(-pz_0)$ over~$C$ is divisible by~$b$. Denote by 
$$
\mathbf{b}: \rPic^{N+g-1}(C) \to \rPic^{d-p}(C)
$$
the map of multiplication by~$b$. We can pull back the line bundle $L(-pz_0)$ to $C \times \rPic^{N+g-1}(C)$ by the map~$\mathbf{b}$ and then extract a $b$th tensor root~$Q$. In other words, we choose a universal line bundle $Q$ over $C \times \rPic^{N+g-1}(C)$ such that 
$$
Q^{\otimes b} = \mathbf{b}^* L(-pz_0).
$$
Thus $Q$ is, just as $L$, a universal line bundle over $C \times \rPic^{d'}(C)$ for $d' = N+g-1$, and the first Chern class of its restriction to $\mbox{point} \times \rPic^{d'}(C)$ vanishes. Therefore we can define, just as for~$L$, the line bundle
$$
\bQ(-\Delta) = \pi_1^*Q \otimes \dots \otimes \pi_N^*Q \left(-\!\!\sum_{1 \leq i<j \leq N}\!\!\Delta_{ij}\right).
$$
over $C^N \times \rPic^{d'}(C)$, that turns out to be a pull-back from~$S^NC \times \rPic^{d'}(C)$. 

Denote by $X_i \in C^N$ the locus where $z_i = z_0$. Then the divisor $X_1 + \dots + X_N$ is the full pullback of the divisor $\xi$ in $S^NC$.

Over $C^N \times \rPic^{d'}(C)$ we have
\begin{align*}
s^*\L(-b \Delta) &=  \pi_1^*s^*L \otimes \dots \otimes \pi_N^*s^*L (-b\Delta) \\
&= \pi_1^*Q^{\otimes b}(p X_1) \otimes \dots \otimes \pi_N^* Q^{\otimes b} (p X_N) (-b \Delta)\\
& = \bQ(-\Delta)^{\otimes b} (p X_1 + \dots + p X_N).
\end{align*}
Thus, on $S^NC \times \rPic^{d'}(C)$ we have
$$
s^*\L(-b\Delta) = \bQ(-\Delta)^{\otimes b}(p \xi).
$$

It follows from the case $b=1$, $p=0$ applied to~$Q$ that 
$$
c_1(\bQ(-\Delta)) = \theta_N - \eta_{N,d'}.
$$
Therefore
$$
c_1(s^*\L(-b\Delta)) = b \theta_N - b \eta_{N,d'} + p \xi
$$
in $H^2(S^NC \times \rPic^{d'}, \Z)$.

Finally, by Eq~\eqref{Eq:pullbacketa} and recalling the identification between $\rPic^d(C)$ and $\rPic^{d-p}(C)$ by $p z_0$ that we used in the beginning of the proof, we get
$$
c_1(\L(-b\Delta)) = b \theta_N - \eta_{N,d} + p \xi.
$$
in $H^2(S^NC \times \rPic^d, \Z)$.
\qed

\section{Applying the GRR formula} \label{Sec:GRR}

Our goal is to apply the Grothendieck-Riemann-Roch formula to the line bundle $\L_b$ and the map
$$
\pi: S^NC \times \rPic^d(C) \to \rPic^d(C).
$$
The formula reads
$$
\ch(R^*\pi_* \L_b) = \pi_*\left(e^{c_1(\L_b)} \, \Td (T_*S^NC) \right).
$$

Concerning the left-hand side, we will show that in the alternating sum
$$
R^*\pi_*\L_b = R^0\pi_*\L_b - R^1\pi_*\L_b + \dots
$$
all the terms vanish except the first one, which is equal to the Laughlin vector bundle. Thus the formula will give us the Chern character of the Laughlin bundle. 

In the right-hand side, $c_1(\L_b) = b \theta_N - \eta + p \xi$ as computed in the previous section. Now we need to compute the Todd class $\Td(T_*S^NC)$ of the tangent vector bundle to $S^NC$.

\subsection{The Kodaira vanishing} \label{Ssec:Kodaira}

In this section we prove that $R^i\pi_* \L_b =0$ for $i \geq 1$. For that it is enough to show that $H^i(S^NC, \L_b) = 0$ for $i \geq 1$ on every fiber $S^NC$ of~$p$. 

\begin{proposition}
The divisor $\theta_N$ in $S^NC$ is nef.
\end{proposition}

\paragraph{Proof.} The divisor $\theta_N$ is known to be ample in the Jacobian of~$C$, therefore it is nef, and so is its pull-back to $S^NC$. \qed

\begin{proposition}
The divisor $\xi$ in $S^NC$ is ample.
\end{proposition}

\paragraph{Proof.} The divisor $z_0$ on $C$ is, of course, ample. It follows that the divisor $X_1 + \dots + X_N$ on $C^N$ (consisting of all $n$-tuples of points where at least one point equals $z_0$) is also ample. (In general, on a product, a sum of ample divisors pulled back from factors is ample. In our case, for $k$ large enough the sections of a $\cO\Bigl(k(X_1 + \dots + X_N)\Bigr)$ that are obtained as products of sections of $\cO(kx_i)$ for each $x_i$ determine an embedding of $C^N$ into $(\CP^{k-g})^N$, which is itself embedded into $\CP^{N(k-g+1)-1}$.)

From this we deduce that $\xi$ is ample on $S^NC$ using the Nakai-Moishezon criterion (\cite{Lazarsfeld}, Thm 1.2.23 and Corollary 1.2.28): a divisor $D$ in $X$ is ample iff for any subvariety $Y \subset X$ of dimension $r$ we have $Y \cdot  D^r >0$. Now, denote by $q: C^N \to S^NC$ the quotient map. Let $Y \subset S^NC$ be a subvariety of dimension~$r$ and $W \subset C^N$ a subvariety that maps finitely onto~$V$ with degree~$d$. Then $Y \cdot \xi^r = \frac1{d} W \cdot (X_1 + \dots + X_N)^r > 0$, because $X_1 + \dots + X_N$ is ample. \qed

\begin{corollary} \label{Cor:ample}
Any divisor $a \xi + b \theta_N$ in $S^NC$ is ample for $a>0$, $b \geq 0$. 
\end{corollary}

\paragraph{Proof.} By Kleiman's criterion (\cite{Lazarsfeld}, Thm 1.4.24) a divisor is ample iff it belongs to the interior of the nef cone. In our case $\xi$, being ample, lies in the interior of the nef cone, while $\theta_N$ lies on the boundary, since it is ample on $\rPic^N(C)$, but has a vanishing intersection number with the projective line in a fiber of $S^NC \to \rPic^N(C)$. Thus $a \xi + b \theta_N$ lies in the interior of the cone for $a>0$, $b \geq 0$. \qed

\begin{proposition} \label{Prop:Kodaira}
If $p >-r$ then the line bundle $\L_b$ on $S^NC$ has no higher cohomology, that is, $h^i(S^NC, \L_b) = 0$ for $i \geq 1$.
\end{proposition}

\paragraph{Proof.} We use the Kodaira vanishing: higher cohomologies of a line bundle $\cO(D)$ vanish if the difference between $D$ and the canonical divisor is ample. We have $c_1(\L_b) = b\theta_N + p \xi$. By Proposition~\ref{Prop:anticanonical}, the canonical class of $S^NC$ is $\theta_N - r \xi$.  The difference equals $(b-1) \theta_N + (p+r) \xi$. By Corollary~\ref{Cor:ample}, this is ample for $b \geq 1$, $p > -r$. \qed

As a corollary, we prove Theorem~\ref{thm:Kodaira}: 
If $p< -r$, then $R^0\pi_*\L_b$ is a vector bundle (the Laughlin bundle), while
$$
R^i \pi_* \L_b = 0 \qquad \mbox{for} \quad i \geq 1.
$$

\paragraph{Proof of Theorem~\ref{thm:Kodaira}.} Given that the sheaf of sections of $\L_b$ is flat over $\rPic^d(C)$, we can use Grauert's theorem (\cite{Hartshorne}, Cor.~12.9). By this theorem, the fact that the spaces $H^i(S^NC, \L_b)$ have the same dimension for any~$L \in \rPic^d(C)$ implies that $R^i \pi_*\L_b$ is a free sheaf over $\rPic^d(C)$ and the spaces $H^i(S^NC, \L_b)$ are its fibers. In other words, $R^0\pi_*\L_b$ is a vector bundle and $R^i \pi_*\L_b$ for $i \geq 1$ is a rank~0 vector bundle. \qed

\begin{corollary} \label{Cor:GRR}
For $p > -r$, the total Chern character of the Laughlin bundle over $\rPic^d(C)$ equals
$$
\pi_* \left[\exp \left(b \theta_N - \eta + p \xi\right) \;  \td^r\xi \; \exp\left(\theta_N\, \frac{\td\,\xi -\xi-1}{\xi}\right) \right].
$$
\end{corollary}

\paragraph{Proof.} This is the Grothendieck-Riemann-Roch formula applied to $\L_b$, given than 
$$
\ch(\L_b) = e^{b\theta_N - \eta + p \xi}
$$
and taking into account that in the alternating sum
$$
\oplus_{i \geq 0} (-1)^i R^i\pi_*\L_b
$$
the first term is the Laughlin bundle and the other terms vanish.
\qed

\subsection{The Todd class of $T_*S^NC$} \label{Ssec:Td}

Recall that $\pi : S^NC \to \rPic^N(C)$ is the projectivization of a vector bundle~$E$ or rank
\begin{equation} \label{Eq:r}
r = N-g+1
\end{equation}
with $c_1(\cO(1)) = \xi$ and $c(E) = e^{-\theta_N}$.

We will use the notation
$$
\td x = \frac{x}{1-e^{-x}}.
$$

\begin{proposition} \label{Prop:anticanonical}
The tangent vector bundle to $S^NC$ has the following characteristic classes:
$$
\begin{array}{rcl}
\ch (T_*S^NC) &=& (g-1)+(r - \theta_N) e^\xi,\\
c(T_*S^NC) & =& (1+\xi)^r \exp \left(-\frac{\theta_N}{1+\xi}\right),\\
\Td (T_*S^NC) &=& (\td \, \xi)^r \exp\left(\theta_N \, 
\frac{\td \, \xi - 1 - \xi}\xi
\right).
\end{array}
$$
\end{proposition}

\paragraph{Proof.} 
Let 
$$
\begin{array}{rcl}
\ln \td \, x &=& \sum\limits_{i \geq 1} a_i x^i,\\
\ln (1+x) &=& \sum\limits_{i \geq 1} b_i x^i,
\end{array}
$$
where, of course, $b_i = (-1)^{i-1}/i$. For any vector bundle~$V$ we have
$$
\begin{array}{rcl}
\ln \Td(V) &=& \sum\limits_{i \geq 1} i! a_i \ch_i(V),\\
\ln \, c(V) &=& \sum\limits_{i \geq 1} i! b_i \ch_i(V).
\end{array}
$$
In particular, we conclude that $\ch_1(E) = -\theta_N$ and $\ch_i(E) = 0$ for $i \geq 2$, because $\ln c(E) = -\theta_N$, so that 
$$
\ch(E) = r - \theta_N.
$$

The tangent vector bundle to $\rPic^N(C)$ is trivial, so it is enough to find the characteristic classes of the relative tangent vector bundle of the family $\pi: S^NC \to \rPic^N(C)$, which is the projectivization of the vector bundle~$E$. The relative tangent bundle is isomorphic to
$$
\Hom(\cO(-1), E/\cO(-1)) = E \otimes \cO(1)/\cO.
$$
The Chern character of $E \otimes \cO(1)$ equals
$$
\ch(E \otimes \cO(1)) = \ch(E) \ch(\cO(1)) = (r-\theta_N) e^\xi,
$$
and the total Chern character of $T_*S^NC$ differs only in the $\ch_0$ term. From the expression for $\ch(E \otimes \cO(1))$ it follows that
$$
c(T_*S^NC) = c(E \otimes \cO(1)) = \exp \left[ \sum_{i \geq 1} i! b_i \ch_i(E \otimes \cO(1)) \right]
$$
$$
=
\exp \left[ r \sum_{i \geq 1} b_i \xi^i - \theta_N \sum_{i \geq 1} b_{i-1} i \xi^{i-1} \right] = \exp(r\ln(1+\xi)) \cdot \exp\left(-\theta_N \frac{d \ln(1+\xi)}{d \xi}\right) 
$$
$$
= (1+\xi)^r \cdot \exp\left(-\frac{\theta_N}{1+\xi}\right).
$$
Similarly, 
$$
\Td(T_*S^NC) = \Td(E \otimes \cO(1)) = \exp \left[ \sum_{i \geq 1} i! a_i \ch_i(E \otimes \cO(1)) \right]
$$
$$
=
\exp \left[ r \sum_{i \geq 1} a_i \xi^i - \theta_N \sum_{i \geq 1} a_{i-1} i \xi^{i-1} \right] = \exp(r\ln(\td \xi)) \cdot \exp\left(-\theta_N \frac{d \ln(\td\, \xi)}{d \xi}\right) 
$$
$$
= \td^r\xi \cdot \exp\left(\theta_N\, \frac{\td\,\xi -\xi-1}{\xi}\right).
$$
\qed

We are now in a position to compute the Chern characters of the Laughlin bundle explicitly: indeed, we have expressed them as intersection numbers on $S^NC \times \rPic^d(C)$ and we have learned enough about the cohomology of this space to compute them. It is worth noting that the computation still looks rather formidable at first sight and that the final answer turns out to be much simpler than what could have been expected. To give the reader a feel of our luck let's consider a concrete simple example.

\begin{example} \label{Ex:luck}
Take $g=0$, $p=0$, $b=1$, $r=6$. Since we are in genus~0, each of the spaces $\rPic^N(C)$ and $\rPic^d(C)$ is reduced to a point, and the classes $\theta_N$ and $\eta$ vanish. The Laughlin bundle is a vector bundle over a point, that is, a vector space. The Chern character equals the dimension of this space. 

The space $S^NC$ is a projective space of dimension $r-1$. The push-forward
$$
\pi_* \left(\td^r(\xi)\right)
$$
(all that remains of our formula!) is the integral over this projective space. The integral of $\xi^{r-1}$ equals 1, the integral of other powers of $\xi$ vanishes for dimensional reasons. Thus the push-forward is equal to the coefficient of $x^{r-1}$ in the power series $\td^r(x)$. 

For $r=6$ we have
$$
\td^6(x) = 
$$
$$
1+3 x +\frac{17}{4} x^2+\frac{15}{4} x^3+
\frac{137}{60} x^4
+ x^5
+\frac{19087}{60480}x^6+
\frac{275}{4032}x^7+\frac{9829}{1209600}x^8-\frac{19}{80640}x^9 + \dots
$$

The coefficient we need is that of $x^5$, and it is strikingly simpler than the other coefficients. To understand the reason of this, let's write the coefficient as a residue:
$$
[x^{r-1}] \td^r(x) = \Res_{x=0} \, (\td(x)/x)^r dx = 
\Res_{x=0} \frac{dx}{ \left(1-e^{-x}\right)^r}.
$$
Now we can perform a change of variables $z = 1 - e^{-x}$, $dz = e^{-x} dx = (1-z) dx$, so that $dz = dz/(1-z)$. The point $x=0$ corresponds to $z=0$. We get
$$
[x^{r-1}] \td^r(x) = 
\Res_{z=0} \frac{dz}{ z^r (1-z)} = [z^{r-1}] \frac1{1-z} = 1.
$$
\end{example}

Drawing hope from this example let's dive into the general case.

\subsection{Computations} \label{Ssec:GRRComputations}

The goal of this section is to provide a closed formula for the expression
$$
\pi_* \left[\exp \left(b \theta_N - \eta + p \xi\right) \;  \td^r\xi \; \exp\left(\theta_N\, \frac{\td\,\xi -\xi-1}{\xi}\right) \right].
$$
from Corollary~\ref{Cor:GRR}. Recall that $\pi$ is the projection $\pi: S^NC \times \rPic^d(C) \to \rPic^d(C)$. Before computing the push-forward of this rather complicated power series, let's compute the pushforward of one monomial.

\begin{lemma} Consider the map $p_2: \rPic^N(C) \times \rPic^d(C) \to \rPic^d(C)$. For $\ell, m \geq 0$, $\ell + m = g$ we have
$$
(p_2)_* (\theta_N^\ell \eta^{2m}) = \frac{\ell! \, (2 m)!}{m!} (-\theta_d)^m.
$$
The push-forward of other monomials ($\ell+ m \not= g$ or odd power of $\eta$) vanish.
\end{lemma}

\paragraph{Proof.} Denote by $(\alpha_i,\beta_i)_{1 \leq i \leq g}$ a symplectic basis of $\rPic^N(C)$ and $(\alpha'_i, \beta'_i)_{1 \leq i \leq g}$ the corresponding symplectic basis of $\rPic^d(C)$. Recall that 
$$
\theta_N = \sum_{i=1}^g \alpha_i \wedge \beta_i, 
\qquad
\eta = \sum_{i=1}^g (\alpha_i \wedge \beta'_i + \alpha'_i \wedge \beta_i),
\qquad
\theta_d = \sum_{i=1}^g \alpha'_i \wedge \beta'_i. 
\qquad
$$
The push-forward selects the terms of the product $\theta_N^\ell \eta^{2m}$ containing, as a factor,
$$
\alpha_1 \wedge \beta_1 \wedge \dots \wedge \alpha_g \wedge \beta_g
$$
and integrates this factor out of the product.

Consider an (anticommuting) monomial $M'$ in variables $\alpha'_i, \beta'_i$ and let's find its coefficient in $(p_2)_* (\theta_N^\ell \eta^{2m})$. Every factor $\alpha'_i$ (respectively, $\beta'_i$) appears in $\eta$ together with $\beta_i$ (respectively, $\alpha_i$). Thus, to a monomial $M'$ corresponds a monomial~$M$ in variables $\alpha_i, \beta$. In order for~$M$ to be completable to 
$$
\alpha_1 \wedge \beta_1 \wedge \dots \wedge \alpha_g \wedge \beta_g
$$
by a power of $\theta_N$, it must contain, for every index~$i$, either both $\alpha_i$ and $\beta_i$ or neither of them. Hence the same should be true of~$M'$. Thus, the only monomials $M'$ that appear in the push-forward with nonzero coefficients have the form 
$$
M'= \alpha'_{i_1} \wedge \beta'_{i_1} \wedge \dots \wedge  \alpha'_{i_m} \wedge \beta'_{i_m}.
$$
In particular, it follows that if the power of $\eta$ is odd, the pushforward vanishes. Now let us find the coefficient of $M'$ in $(p_2)_* (\theta_N^\ell \eta^{2m})$.

Each factor of $M'$ can be selected in any of the $2 m$ factors $\eta$. There are $(2 m)!$ choices like that. We also acquire a sign $(-1)^m$ by permuting the factors of $\eta^{2 m}$ to separate $M$ and $M'$. For instance,
$$
(\alpha_1 \wedge \beta'_1) \wedge (\alpha'_1 \wedge \beta_1) =
- (\alpha_1 \wedge \beta_1) \wedge (\alpha'_1 \wedge \beta'_1).
$$
 (Since both $M$ and $M'$ have even degrees, it does not matter which of them we put first.)

Further, there is a unique monomial $M^{\rm comp}$ that completes $M$ to 
$$
\alpha_1 \wedge \beta_1 \wedge \dots \wedge \alpha_g \wedge \beta_g
$$
and it appears in $\theta_N^\ell$ with coefficient $\ell!$, provided $\ell+m=g$. (Otherwise, the monomials of $\theta_N^\ell$ are of wrong degree, and the push-forward vanishes.) Thus $M'$ appears in the push-forward with coefficient
$$
(-1)^m \, (2 m)! \, \ell!.
$$
Finally, we note that $\theta_d^m$ contains exactly the same monomials $M'$, but with coefficient $m!$. Thus the answer can be more compactly written as
$$
(p_2)_* (\theta_N^\ell \eta^{2m}) = \frac{s! \, (2 m)!}{m!} (-\theta_d)^m
$$
provided $\ell+m=g$.
\qed

\begin{lemma} \label{Lem:monomial}
For $k,\ell,m \geq 0$, $k+\ell+m=g$ we have
$$
\pi_* (\xi^{N-g+k} \theta_N^\ell \eta^{2m}) = \frac{(\ell+k)! \, (2 m)!}{k!\,m!} \; (-\theta_d)^{m}.
$$
The pushforwards of all other monomials ($k+\ell+m \not= g$, odd power of $\eta$ or power of $\xi$ below $N-g$) vanish.
\end{lemma}

\paragraph{Proof.}
Decompose $\pi$ into a composition of two maps
$$
S^NC \times \rPic^d(C) \stackrel{p_1}{\longrightarrow} \rPic^N(C) \times \rPic^d(C) \stackrel{p_2}{\longrightarrow} \rPic^d(C).
$$
First apply the push-forward under $p_1$. We have ${(p_1)}_* \xi^{N-g+k} = W_k = \theta_N^k / k!$ for $k \geq 0$. The push-forward of lower powers of $\xi$ vanishes. As for the classes $\theta_N$ and $\eta$ they are defined as pullbacks under~$p_1$. Thus, by the projection formula,
$$
(p_1)_* (\theta_N^\ell \eta^{2m} \xi^{N-g+k}) = 
\frac1{k!} \theta_N^{\ell+k} \eta^{2m}.
$$
Now apply the previous lemma. \qed

Using this lemma we can isolate the coefficient of $\theta_d^m$ in the Grothendieck-Riemann-Roch formula
$$
\pi_* \left[\exp \left(b \theta_N - \eta + p \xi\right) \;  \td^r\xi \; \exp\left(\theta_N\, \frac{\td\,\xi -\xi-1}{\xi}\right) \right].
$$
Specifically, we write
\begin{equation} \label{Eq:Chernm}
\ch_m(V) =\frac1{(2m)!}\;
\pi_* \left[\eta^{2m} \cdot e^{p \xi} \;  \td^r\xi \cdot \exp\left(\theta_N\,\left[ \frac{\td\,\xi -\xi-1}{\xi}+b \right]\right) \right].
\end{equation}

\begin{lemma} \label{Lem:AB}
Let $A(\xi)$ and $B(\xi)$ be two power series. Then 
$$
\frac1{(2m)!} \pi_* \left[\eta^{2m} \cdot A(\xi) \cdot \exp(\theta_N B(\xi))  \right]= 
\frac{(-\theta_d)^m}{m!} \cdot [x^{N-g}] A(x) \left(B(x) + \frac1x \right)^{g-m},
$$
where $[x^{N-g}]$ is the coefficient of $x^{N-g}$ in the Laurent series.
\end{lemma}

\paragraph{Proof.} Decompose the left-hand side by powers of $\theta_N$ from 0 to $g-m$ (the push-forward of higher powers of $\theta_N$ times $\eta^{2m}$ vanish). We get
$$
\frac1{(2m)!} \, \sum_{\ell=0}^{g-m} \frac1{\ell!} \cdot \pi_* \left( \eta^{2m} \theta_N^\ell \xi^{N-m-\ell} \right)
\cdot [x^{N-m-\ell}] A(x) B^\ell(x). 
$$
Indeed, if we fix the power of $\theta_N$ to be $\ell$, that means that we choose the term $\left(\theta_N B(\xi)\right)^\ell/\ell!$ in the exponential. And $[x^{N-m-\ell}] A(x) B^\ell(x)$ is the coefficient of the appropriate power of $\xi$ in $A(\xi) B(\xi)^\ell$.

Now, by Lemma~\ref{Lem:monomial}, we replace $\pi_* \left( \eta^{2m} \theta_N^\ell \xi^{N-m-\ell} \right)$ with
$$
\frac{(g-m)! \, (2 m)!}{(g-m-\ell)!\,m!} \; (-\theta_d)^{m}.
$$
Simplifying by $(2m)!$ we get
$$
\frac1{(2m)!} \pi_* \left[\eta^{2m} \cdot A(\xi) \cdot \exp(\theta_N B(\xi))  \right]
= 
\frac{(-\theta_d)^m}{m!} \sum_{\ell=0}^{g-m}
{g-m \choose \ell}  [x^{N-m-\ell}] A(x) B^\ell(x).
$$

This can by rewritten as 
$$
\frac{(-\theta_d)^m}{m!} [x^{N-g}] A(x) \left(B(x)+\frac1x\right)^{g-m},
$$
which is the right-hand side of the equality of the lemma. \qed

\begin{corollary} \label{Cor:powerg}
We have the following expression for the Chern characters of the Laughlin bundle:
$$
\ch_m(V)  = \frac{(-\theta_d)^m}{m!} \;
 [x^{r-1}] \, e^{px} \, \td^rx \, \left(\frac{\td \,x}x+b-1\right)^{g-m}.
$$
\end{corollary}

This follows immediately by plugging into Equation~\eqref{Eq:Chernm} the result of Lemma~\ref{Lem:AB} and recalling that $N-g = r-1$. \qed

\begin{lemma} \label{Lem:residue}
We have
$$
[x^{r-1}] \, e^{px} \, \td^rx \, \left(\frac{\td \,x}x+b-1\right)^{g-m}
= \sum_{k=m}^g
{{g-m}\choose{k-m}}{{N-g+p}\choose{k-g+p}}b^{k-m}.
$$
\end{lemma}

\paragraph{Proof.}
To prove the lemma we use the same trick as in Example~\ref{Ex:luck}: represent the coefficient as a residue and perform a change of variables $z = 1- e^{-x}$. Then
$$
\td\,x/x = 1/z, \qquad, e^{px} = \frac1{(1-z)^p}, \qquad dx = \frac{dz}{1-z}.
$$
We get
$$
[x^{r-1}] \, e^{px} \, \td^rx \, \left(\frac{\td \,x}x+b-1\right)^{g-m} = 
\Res_{x=0} \, e^{px} \, \left( \frac{\td\, x}x \right)^r \, \left(\frac{\td \,x}x+b-1\right)^{g-m}\, dx 
$$
$$
= \Res_{z=0} \frac1{(1-z)^p\, z^r} \, (1/z+b-1)^{g-m} \frac{dz}{1-z} = 
[z^{r-1}] \frac{\left(b+ \frac{1-z}z\right)^{g-m}}{(1-z)^{p+1}}.
$$
The last coefficient is easy to compute. We have
$$
\left(b+ \frac{1-z}z\right)^{g-m} = 
\sum_{k=0}^{g-m} {g-m \choose k} b^{g-m-k} \frac{(1-z)^k}{z^k}.
$$
Thus
$$
[z^{r-1}] \frac{\left(b+ \frac{1-z}z\right)^{g-m}}{(1-z)^{p+1}} =
\sum_{k=0}^{g-m} {g-m \choose k} b^{g-m-k} [z^{r-1+k}] \frac1{(1-z)^{p+1-k}}.
$$

Now,
$$
[z^q] \frac1{(1-t)^s} = {q+s-1 \choose s-1}.
$$
Plugging $q=r-1+k= N-g+k$, $s= p+1-k$ we get
$$
[z^{r-1}] \frac{\left(b+ \frac{1-z}z\right)^{g-m}}{(1-z)^{p+1}} =
\sum_{k=0}^{g-m} {g-m \choose k} b^{g-m-k} {N-g+p \choose p-k}.
$$
To put the answer in the final form, we change the summation index from $k$ to $g-k$ and rewrite the last sum as
$$
\sum_{k=m}^g {g-m \choose g-k} {N-g+p \choose k-g+p} b^{k-m}.
$$

\qed

Theorem~\ref{thm:ch} (the expression for Chern characters of the Laughlin bundle) follows from Corollary~\ref{Cor:powerg} and Lemma~\ref{Lem:residue}.

\end{document}